\documentclass[12pt,a4paper,twoside]{article}

\usepackage{a4, amsmath, amssymb, bm}

\newcommand{\eps}{\varepsilon}
\def\NN{\mathbb N}
\def\RR{\mathbb R}
\def\HH{\mathbb H}
\def\calf{{\cal F}}
\def\ric{\operatorname{Ric}}
\def\tr{\operatorname{Tr}}
\newcommand\vol{\operatorname{vol}}
\newcommand{\eq}{\hspace*{-2mm}&=&\hspace*{-2mm}}
\newcommand{\plus}{\hspace*{-2mm}&+&\hspace*{-2mm}}

\newtheorem{df}{Definition}
\newtheorem{cor}{Corollary}
\newtheorem{ex}{Example}

\newtheorem{prop}{Proposition}
\newtheorem{thm}{Theorem}
\newtheorem{rem}{Remark}

\title{On the partial Ricci curvature of foliations}

\author{Vladimir Y. Rovenski
\footnote{The work was supported by grant P-IEF, No. 219696 of Marie-Curie action.
\newline
Author acknowledges Prof. Pawel Walczak (University of Lodz) for help and useful discussion of the work. Theorem~\ref{T-Kmix} (in Section~\ref{sec:scalar}) is the result of cooperation with him.}
\\
         \small Department of Mathematics,
             Faculty of Science and Science Education\\
          \small University of Haifa, Mount Carmel, Haifa, 31905, Israel\\
          \small E-mail: rovenski@math.haifa.ac.il}

\date{}

\begin{document}

\maketitle

\begin{abstract}
We consider a problem of prescribing the partial Ricci curvature on
a locally conformally flat manifold $(M^n, g)$ endowed with the complementary orthogonal distributions ${\mathcal{D}}_1$ and ${\mathcal{D}}_2$. We provide conditions for symmetric $(0,2)$-tensors $T$
of a simple form (defined on $M$) to admit metrics $\tilde g$, conformal to $g$, that solve the partial Ricci equations. The solutions are given explicitly. Using above solutions,
we also give examples to the problem of prescribing the mixed scalar curvature related to ${\mathcal{D}}_i$. In aim to find ''optimally placed" distributions, we calculate the variations of the total mixed scalar curvature (where again the partial Ricci curvature plays a key role), and give examples concerning minimization of a total energy and bending of a distribution.

\noindent
\textbf{Keywords/Phrases}:
Riemannian manifold; distribution; conformal; variation; partial Ricci curvature;
mixed scalar curvature;
energy.

\end{abstract}

\section*{Introduction}

Let ${\mathcal{D}}_{1}$ and ${\mathcal{D}}_{2}$ be smooth complementary orthogonal distributions on a~Riemannian manifold $(M^n, g)$ with the Levi-Civita connection $\nabla$ and the Riemannian curvature tensor $R$.
Assume that $\dim {\mathcal{D}}_i=p_i>0\ (i=1,2)$, $\dim M=n,\ n=p_1+p_2$.
Let~$e_1,\dots e_n$ be a local orthonormal frame adapted to ${\mathcal{D}}_{1}$ and ${\mathcal{D}}_{2}$, i.e., $e_i\in {\mathcal{D}}_{1}$ for $i\le p_1$ and $e_\alpha\in {\mathcal{D}}_{2}$ for $\alpha>p_1$.
The~\textit{mixed scalar curvature} is given by
$K_{1,2}=\sum\nolimits_{i,\alpha} g(R(e_i,e_\alpha)e_\alpha,e_i)$, see~\cite{w90}.
We~call the \textit{partial Ricci curvature related to ${\mathcal{D}}_1$} the symmetric bilinear form
 $\ric_{1}(X,Y)=\sum\nolimits_{\,\alpha>p_1}g(R(e_\alpha, X)Y, e_\alpha)$
on the tangent bundle $TM$. The definition for ${\mathcal{D}}_{2}$ are similar.
Indeed, $K_{1,2}=\tr_g\ric_{i|\,{\mathcal{D}}_i}$ for $i=1,2$.
We are interested in the problem that concerns ``optimally placed" distributions

\vskip1mm
($\bm P_0$) \textit{Find variational formulae for the functional}
 $I_K: {\mathcal{D}}_1\to\int_M K_{1,2}\,{\rm d}\vol$.
\vskip1mm

\noindent
The variational formulae for $I_K$ (and for $i$-th mean curvatures) related to  codimension one distributions on a compact $(M, g)$ are developed in \cite{rw4}.
In the paper we consider distribution ${\mathcal{D}}_1$ of any codimension,
represent the first and second variations of the total mixed scalar curvature
and give examples concerning minimization of the total energy and bending of distributions.
The partial Ricci curvatures
play a key role in above variational formulae.

 Different aspects of the problem of finding a metric $g$, whose Ricci tensor is a given second-order symmetric tensor $T$, were considered by several authors, see \cite{bes}, \cite{pt06}\,--\,\cite{pt09}, etc.
The following problems for a differentiable manifold $M$ with transversal complementary distributions ${\mathcal{D}}_{1}$ and ${\mathcal{D}}_{2}$, which generalize the classical ones, seem to be interesting:

\vskip1mm
($\bm P_1$) \textit{Given a symmetric $(0,2)$-tensor $T$ on $M$ satisfying $T({\mathcal{D}}_1,{\mathcal{D}}_2)=0$, does there exist a Riemannian metric $g$ with the properties $g({\mathcal{D}}_{1},{\mathcal{D}}_{2})=0$ and
either $\ric_{i\,|\,{\mathcal{D}}_{i}}= T_{|\,{\mathcal{D}}_{i}}$
or $[\ric_{i}-\frac12 K_{1,2}\,g]_{\,|\,{\mathcal{D}}_{i}}= T_{|\,{\mathcal{D}}_{i}}$, where $i=1,2$\,?}

\vskip1mm
($\bm P_2$) \textit{Given a function $\bar K\in C(M)$,
does there exist a Riemannian metric $g$ on $M$, whose mixed scalar curvature
(related to ${\mathcal{D}}_{1}$ and ${\mathcal{D}}_{2}$) is $\bar K$?}

\vskip1mm
Note that ($\bm P_1$) for $T=0$ asks about existence of either ${\mathcal{D}}_i$-flat or ``${\mathcal{D}}_i$-Einstein" metrics.
($\bm P_2$)~is similar to the known problem of prescribing scalar curvature on $M$;
its particular case $\bar K=const$ corresponds to the \textit{Yamabe problem} (of prescribing constant scalar curvature on~$M$).

We study the problems ($\bm P_1$) and ($\bm P_2$) on a locally conformally flat $(M,g)$,
in particular, on space forms, for tensors $T$ of a simple form.
We find necessary and sufficient conditions on $T$ for the existence of metrics $\tilde g=(1/{\phi\,^2})\,g$ (conformal to the metric $g$) which solve the systems
\begin{equation}\label{E-PDE-c1}
 a)\ \widetilde\ric_{i\,|\,{\mathcal{D}}_{i}}= T_{|\,{\mathcal{D}}_{i}}\ (i=1,2);\quad
 b)\ [\,\widetilde\ric_{i}-\frac12\widetilde K_{1,2}\,\tilde g]_{\,|\,{\mathcal{D}}_{i}}= T_{|\,{\mathcal{D}}_{i}}\ (i=1,2).
\end{equation}
The compatibility condition for (\ref{E-PDE-c1})(a) is $\tr_g T_{|\,{\mathcal{D}}_{1}}=\tr_g T_{|\,{\mathcal{D}}_{2}}$ with the traces equal to $K_{1,2}$, while for (\ref{E-PDE-c1})(b) is
$(1- p_2/2)\tr_g T_{|\,{\mathcal{D}}_{1}}=(1- p_1/2)\tr_g T_{|\,{\mathcal{D}}_{2}}$.

\vskip1mm
In \textbf{Section~\ref{sec:ricci}} we determine all tensors $T$ of ($\bm P_1$), the functions $\bar K$ of ($\bm P_2$) and the corresponding metrics $\tilde g$ that solve the systems (\ref{E-PDE-c1}). Theorems~\ref{T-main-1}\,--\,\ref{T-main-4} and Corollaries~\ref{Cor-main1}\,--\,\ref{C-main1} extend recent results of \cite{pt06}\,--\,\cite{pt09} (where ${\mathcal{D}}_1=TM$ and ${\mathcal{D}}_2=0$) to cases of the partial Ricci
and the mixed scalar curvatures of distributions.
 In~\textbf{Section~\ref{sec:scalar}}, in aim to find ''optimally placed" distributions, see ($\bm P_0$), we calculate the first and second variations of total $K_{1,2}$ using the partial Ricci curvature, and give examples concerning minimization of a total energy and bending of a distribution.
 \textbf{Section~\ref{sec:proof}} contains proofs of results.

\section{Prescribed partial Ricci curvature}
\label{sec:ricci}

We start with the solution to (\ref{E-PDE-c1}), see ($\bm P_1$), at a point $x\in M$.
(The~constant curvature metrics are solutions to
$\ric_{i\,|\,{\mathcal{D}}_{i}}=\lambda_i\,g_{\,|\,{\mathcal{D}}_{i}}$ at one point).
Let $M$ be a neighborhood of the origin $O$ in $\RR^n=\RR^{p_1}\times\RR^{p_2}$,
and ${\mathcal{D}}_{i}\ (i=1,2)$ is tangent to the $i$-th factor.

\begin{prop}\label{P-main1}
 Let $T$ be a diagonal $n$-by-$n$ matrix satisfying condition
\[
 K=\sum\nolimits_{i\le p_1} T_{ii}=\sum\nolimits_{\alpha>p_1} T_{\alpha\alpha}.
\]
 Then the metric $g$ in a neighborhood of $O$ in $\RR^n=\RR^{p_1}\times\RR^{p_2}$
 satisfying (\ref{E-PDE-c1})(a) at the origin, can be selected in the form
 $g = \sum_{\,a=1}^{n}\big(1-\sum_{\,b\le n} c_{bb} x_b^2\big) dx^a\otimes dx^a$,
 where, for example,
 \[
 c_{ii}=T_{ii}/p_2\ \ (1\le i\le p_1),\quad
 c_{\alpha\alpha}=(T_{\alpha\alpha}-K/p_2)/p_1\ \ (p_1<\alpha\le n).
 \]
\end{prop}

Next we will represent the formulae relating partial Ricci and mixed scalar curvatures for two conformally related metrics on $M$ with distributions ${\mathcal{D}}_1,{\mathcal{D}}_2$.
(Notice that the conformal change of a metric $g$ preserves the orthogonality of ${\mathcal{D}}_{1}$ and ${\mathcal{D}}_{2}$.)

We~call $\Delta^{(1)}\phi=\sum_{i\le p_1} h_\phi(e_i, e_i)$
and $\Delta^{(2)}\phi=\sum_{\alpha>p_1} h_\phi(e_\alpha, e_\alpha)$ the ${\mathcal{D}}_{1}$- and ${\mathcal{D}}_{2}$-\textit{laplacian} of $\phi$, respectively, where $h_\phi$ is the Hessian of $\phi$.
The Hessian of $\phi$ is the symmetric $(0,2)$-tensor $h_\phi(X,Y)=g(S(X),Y)$,
where $S(X)=\nabla_X\nabla\phi$ is a~self-adjoint $(1,1)$-tensor, and $\nabla\phi$ is the gradient of $\phi$.
Indeed, $\Delta^{(1)}+\Delta^{(2)}=\Delta$.
 Let $\RR^n=\RR^{p_1}\times\RR^{p_2}\ (p_1,p_2>0)$ be a decomposition of Euclidean space.
The partial laplacians are
$\Delta^{(1)}=\sum\nolimits_{\,i\le p_1} \frac{\partial^2}{\partial x_i^2}$ and
$\Delta^{(2)}=\sum\nolimits_{\,\alpha>p_1}\frac{\partial^2}{\partial x_\alpha^2}$.

\begin{prop}\label{E-RicD}
Let $(M,g)$ be a Riemannian manifold with complementary orthogonal distributions ${\mathcal{D}}_{1},{\mathcal{D}}_{2}$, and $\phi:M\to\RR_+$ a smooth function.
The partial Ricci curvatures and the mixed scalar curvature transform under conformal change of a metric $\tilde g=(1/\phi^2)g$ by the formulae
 \begin{eqnarray}\label{E-Rij3a}
 \nonumber
 \widetilde\ric_{1} \eq \ric_{1}
 +\big[p_2\phi\,h_{\phi}+\big(\phi\,\Delta^{(2)}\phi-p_2\,|\nabla\phi|^2\big)g\big]/{\phi\,^2},\\
 \widetilde\ric_{2} \eq \ric_{2}+
 \big[p_1\phi\,h_{\phi}+\big(\phi\,\Delta^{(1)}\phi-p_1\,|\nabla\phi|^2\big)g\big]/{\phi\,^2},\\
\label{E-tildeK1}
 \widetilde K_{1,2} \eq \phi^2 K_{1,2}
 +\phi\,(p_1\,\Delta^{(2)}\phi +p_2\,\Delta^{(1)}\phi)-p_1p_2|\nabla\phi|^2.
\end{eqnarray}
\end{prop}

\begin{cor}\label{Cor-T1-u}
Given complementary orthogonal distributions ${\mathcal{D}}_{1}$ and ${\mathcal{D}}_{2}$
of equal dimensions $p_1=p_2>0$  on $(M,g)$ with mixed scalar curvature $K$, let $\bar g=u^{4/(n-2)}g$ be the conformal metric, where $u>0$ is a function on $M$.
Then the mixed scalar curvature $\bar K$ of $\,\bar g$
satisfies the PDE
\begin{equation}\label{E-PDE-RnKK}
 -\frac{n}{n-2}\,\Delta u +K u=\bar K u^{\frac{n+2}{n-2}}.
\end{equation}
\end{cor}

\begin{rem}\rm One may extend above formulae for a pseudo-Riemannian metric.
The formulae (\ref{E-Rij3a}) and (\ref{E-PDE-RnKK}) are similar to the classical formulae
for the Ricci and scalar curvatures (see, for example, \cite{pt07})
 \begin{eqnarray*}
 \widetilde\ric -\ric = \big[(n-2)\phi\,h_{\phi}+\big(\phi\,\Delta\phi -(n-1)\,|\nabla\phi|^2\big)g\big]/{\phi\,^2},\\
 -4\,\frac{n-1}{n-2}\,\Delta u +K u=\bar K u^{\frac{n+2}{n-2}},\quad
 \bar g=u^{4/(n-2)} g.
\end{eqnarray*}
\end{rem}

We will consider the problem ($\bm P_1$) for a neighborhood $V\subset M^n$
of a locally conformally flat space. Suppose that there are coordinates $(x_1,\ldots,x_n)$ on $V$ with the metric $g_{i j}=\delta_{i j}/F^2$, where $F>0$ is a differentiable function on $M$. We will fix on $V$ canonical foliations
${\cal F}_1=\{x_\alpha=c_\alpha,\ \alpha>p_1\}$ and
${\cal F}_2=\{x_i=c_i,\ i\le p_1\}$, $c_i,c_\alpha\in\RR$, consisting of coordinate submanifolds.
Let ${\mathcal{D}}_{2}=T{\cal F}_1$ and ${\mathcal{D}}_{1}=T{\cal F}_2$ be their tangent distributions.

\begin{thm}\label{T-main-1}
Let $(M^n,\bar g)$ be a locally conformally flat Riemannian manifold with complementary orthogonal distributions ${\mathcal{D}}_{1},{\mathcal{D}}_{2}$ of dimensions $p_1,p_2\ge2$,
and $V\subset M^n$ an open set with coordinates $(x_1,\ldots,x_n)$ such that
$\bar g_{ij}=\delta_{i j}/F^2$.
Suppose that $T$ is a symmetric $(0,2)$-tensor with the properties
\begin{equation}\label{E-prop-th1}
 T_{ij}=f_1\delta_{ij},\quad
 T_{\alpha\beta}=f_2\delta_{\alpha\beta},\quad
 T_{i\alpha}=0
 \quad(i,j\le p_1,\ \alpha,\beta>p_1),
\end{equation}
where $f_1,f_2\in C^1(V)$. Then, in any of cases (a) or (b), there is a metric $\tilde g=(1/\phi^2)\bar g$ solving the problem (\ref{E-PDE-c1}) if and only if
 $\phi\,F=\sum\nolimits_{i\le p_1}(a_1 x_i^2+b_i x_i)
 +\sum\nolimits_{\alpha>p_1}(a_2 x_\alpha^2+b_\alpha x_\alpha)+c$,
 and
\begin{eqnarray}\label{E-Th1-a}
 a) &&\hskip-5mm
 f_1=-\frac{p_2}{(\phi\,F)^2}[\lambda-2(a_2-a_1)\mu],\quad
 f_2=-\frac{p_1}{(\phi\,F)^2}[\lambda-2(a_2-a_1)\mu],\\
 \label{E-Th1-b}
 b) &&\hskip-5mm
 f_1=\frac{p_2(p_1\!-2)}{2\,(\phi\,F)^2}[\lambda-2(a_2\!-a_1)\mu],\ \
 f_2=\frac{p_2(p_1\!-2)}{2\,(\phi\,F)^2}[\lambda-2(a_2\!-a_1)\mu].\quad
\end{eqnarray}
 Here $a_1, a_2, b_k, c\in\RR$, and
\begin{equation}\label{E-lambda-mu}
 \lambda=(\sum\nolimits_{k} b_k^2)-2(a_1+a_2)\,c,\ \
 \mu=\sum\nolimits_{i\le p_1}(a_1 x_i^2+b_i x_i)
 +\!\sum\nolimits_{\alpha>p_1}(a_2 x_\alpha^2+b_\alpha x_\alpha).
\end{equation}
\end{thm}

\begin{cor}\label{Cor-main1}
Let ${\mathcal{D}}_{1},{\mathcal{D}}_{2}$ be tangent distributions to canonical foliations on the Euclidean product space $(\RR^{n}=\RR^{p_1}\times\RR^{p_2},g)$, and $p_1,p_2\ge 2$. Suppose that $T$ is a symmetric $(0,2)$-tensor satisfying (\ref{E-prop-th1}),
where $f_1,f_2\in C^1(\RR^{n})$. Then there is a metric $\tilde g=(1/\phi^2)g$ solving the problem (\ref{E-PDE-c1}) (in any of cases (a) or (b)) if and only~if
 $\phi=\sum\nolimits_{i\le p_1}(a_1 x_i^2+b_i x_i)
 +\sum\nolimits_{\alpha>p_1}(a_2 x_\alpha^2+b_\alpha x_\alpha)+c$, and
 $f_a,\,f_\beta$ are given by  (\ref{E-Th1-a}), (\ref{E-Th1-b}),
 where $a_1, a_2, b_k, c\in\RR$, and $\lambda,\,\mu$ are defined in (\ref{E-lambda-mu}).
A~non-complete metric $\tilde g$ is defined on $\RR^{n}$,
if either $a_1,a_2>0$ and $\frac1{a_1}\sum_i b_i^2+\frac1{a_2}\sum_\alpha b_\alpha^2<4c$
or $a_1,a_2<0$ and $\frac1{a_1}\sum_i b_i^2+\frac1{a_2}\sum_\alpha b_\alpha^2>4c$.
In other cases, excluding the homothety, $\tilde g$ has singular points.
\end{cor}

\begin{ex}\label{R-a1-a2}\rm
If $a_i=a$ in Corollary \ref{Cor-main1}, then $\widetilde K_{1,2}=-p_1p_2\lambda$ defined in (\ref{E-lambda-mu}), and the singularity set of $\tilde g$ can be explicitly described in terms of~$\lambda$. Namely,
if $\lambda<0$ then a non-complete metric $\tilde g$ is defined on $\RR^{n}$, and
if $\lambda\ge0$ then, excluding the homothety, the set of singularity points of $\tilde g$ consists of
 $1)$~a~point if $\lambda=0$;
 $2)$~a~hyperplane if $\lambda>0$ and $a=0$;
 $3)$~an~$(n-1)$-dimensional sphere if $\lambda>0$ and $a\ne0$.

  One may consider the pseudo-Euclidean space $(\RR^n, g)$ with the coordinates $x=(x_1,\ldots,x_n)$ and the metric $g_{km}=\epsilon_k\delta_{km},\ \epsilon_k=\pm1$.
  Then, for example, in case (a) of Corollary~\ref{Cor-main1} for $a_1=a_2=a$, we have $\phi=\sum\nolimits_{k=1}^{n}(\epsilon_k a x_k^2+b_k x_k)+c$,
where $p_1\sum\nolimits_{\alpha>p_1}\epsilon_\alpha=p_2\sum\nolimits_{i\le p_1}\epsilon_i$,
and
 $f_1=-p_2\lambda/\phi^2-(2a/\phi)(p_2-\sum\nolimits_\alpha\epsilon_\alpha)$,
 \
 $f_2=-p_1\lambda/\phi^2-(2a/\phi)(p_1-\sum\nolimits_i\epsilon_i)$.
\end{ex}

\begin{ex}\rm
We will discuss our results, when $M$ is the hyperbolic space $(\HH^{n}, \bar g)$, represented by the half space model $\RR^n_+,\ x_n>0$, and $\bar g_{ij}=\delta_{ij}/x_n^2$. For any pair of integers $p_1,p_2>0$, $p_1+p_2=n$, denote by $\calf$ a foliation by $p_2$-planes $\{x\}\times\RR^{p_2}$, where $x=(x_1,\ldots,x_{p_1},0,\ldots,0)$. Let ${\mathcal{D}}_{2}$ be the distribution tangent to $\calf$, and ${\mathcal{D}}_{1}$ its orthogonal complement.
 Using $F=x_n$ in Theorem~{\ref{T-main-1}}, we obtain
$\phi\,x_n=\sum\nolimits_{i\le p_1}(a_1 x_i^2+b_i x_i)
 +\sum\nolimits_{\alpha>p_1}(a_2 x_\alpha^2+b_\alpha x_\alpha)+c$.

If $a_1=a_2=a$, the mixed scalar curvature $\widetilde K_{1,2}=-p_1p_2\lambda$. Moreover, in this case a non-complete metric $\tilde g$ is defined on $\HH^{n}$ whenever
 (i) $\lambda<0$;
 (ii) $\lambda=0$ and $a=0$ (hence $c\ne0$);
 (iii)~$\lambda=0,\ a\ne0$ and $b_n/a\ge0$;
 (iv) $\lambda>0,\ a=0$, $b_k=0$ for $k<n$ and $c/b_n\ge0$;
 (v)~$\lambda>0,\ a\ne0$ and $b_n/a\ge\sqrt\lambda/|a|$.
Otherwise, the singularity set of $\tilde g$ consists of intersection of a hyperplane or a sphere with the half-space $x_n>0$.

We show (i)--(v) in case a). From (\ref{E-tildeK1}) and $\ric_{i}=0$, we get $\widetilde K_{1,2}=-p_1p_2\lambda$.

 If $\lambda<0$ then $\tilde g$ is defined on $\HH^n$ and $T$ is positive definite.
 Let $\lambda=0$. If $a=0$ then $b_k=0$ for $k\le n$, $\phi=c/x_n\ne0$ and $\tilde g$ is defined on $\HH^n$. If $a\ne0$, then if $b_n/(2a)\ge0$ then $\tilde g$ is defined on $\HH^n$; otherwise, if $b_n/(2a)<0$ then $\tilde g$ has a singularity at the point
 $\tilde x=-(b_1,\ldots b_n)/(2a)$.

 Let $\lambda>0$. If $a=0$, we have two cases.
 In the first one, we have $b_k=0$ for $k<n$ and $b_n\ne0$.
 In this case, if $c/b_n\ge0$ then $\tilde g$ is defined on $\HH^n$; otherwise if $c/b_n<0$, then any point of the hyperplane $x_n=-c/b_n$ is a singularity point of $\tilde g$.
 In the second case, we have $b_{k_0}\ne0$ for some $k_0<n$, then any point that belongs to the intersection of the hyperplane $(\sum_k b_k x_k)+c=0$ with the half-space $x_n>0$, is a singularity point of $\tilde g$.

When $\lambda>0$ and $a\ne0$, if $b_n/a\ge\sqrt\lambda/|a|$, then $\tilde g$ is defined on $\HH^n$. Otherwise, if $b_n/a<\sqrt\lambda/|a|$, then any point $p$ of the $(n-1)$-dimensional sphere centered at the point of coordinates
$\tilde x=-(b_1,\ldots b_n)/(2a)$, with radius $\sqrt\lambda/(2|a|)$, such that $p$ is the half space $x_n>0$, is a point of singularity of $\tilde g$.
The case b) is similar to case a).
\end{ex}

The next theorem extends Theorem~\ref{T-main-1} to the tensors $T$ such that
$T_{|{\mathcal{D}}_1}=\sum\nolimits_{\,i\le p_1} f_i(x_k)\,dx_i^2$ and
$T_{|{\mathcal{D}}_2}=\sum\nolimits_{\,\alpha>p_1} f_\alpha(x_k)\,dx_\alpha^2$, with a fixed index $k\le p_1$.

\begin{thm}\label{T-main-2}
Let $(M^n,\bar g)$ be a locally conformally flat Riemannian manifold with complementary orthogonal distributions ${\mathcal{D}}_{1},{\mathcal{D}}_{2}$ of dimensions $p_1,p_2\ge3$.
Let $V\subset M^n$ be an open set with coordinates $(x_1,\ldots,x_n)$ such that
$\bar g_{ij}=\delta_{i j}/F^2$. Consider a symmetric $(0,2)$-tensor $T$ satisfying
\begin{equation}\label{E-prop-th2}
 T_{ij}=f_i(x_k)\,\delta_{ij},\quad
 T_{\alpha\beta}=f_\alpha(x_k)\,\delta_{\alpha\beta},\quad
 T_{i\alpha}=0
 \quad(i,j\le p_1,\ \alpha,\beta>p_1),
\end{equation}
with a fixed $k\le p_1$.
Suppose that the functions $f_i,f_\alpha\in C^1(V)$
and $f_i$ are not all equal.
Then, in any of cases (a) or (b), there is a metric $\tilde g=\bar g/\phi^2$ solving the problem (\ref{E-PDE-c1}) if and only~if there is a differentiable function $U(x_k)$ on $V$ such that $\phi F=e^{U}$ and
\begin{eqnarray}\label{E-T4-fh-a}
\nonumber
 a) &&\hskip-5mm f_k=p_2 U'',\ \ f_i=-p_2{U'}^2\ \ (i\ne k),\ \
       f_\alpha=U''-(p_1-1){U'}^2\ \ (\alpha > p_1),\\
\nonumber
 b) &&\hskip-5mm f_k=\frac12\,p_2\big(U''+(p_1-1){U'}^2\big),\ \
       f_i=-\frac12\,p_2\big(U''-(p_1-3){U'}^2\big)\ \ (i\ne k),\\
    &&\hskip-5mm f_\alpha=\frac12\,(p_2-2)[(p_1-1){U'}^2-U'']\quad (\alpha > p_1).
\end{eqnarray}
\end{thm}

\begin{cor}\label{C-main-33}
Let ${\mathcal{D}}_{1},{\mathcal{D}}_{2}$ be tangent distributions to canonical foliations on the
Euclidean product space $(\RR^{n}=\RR^{p_1}\times\RR^{p_2},\,g)$ with $p_1,p_2\ge 3$.
Consider a symmetric $(0,2)$-tensor $T$ with the properties (\ref{E-prop-th2}) for a fixed $k\le p_1$.
Suppose that the functions $f_i,f_\alpha\in C^1(\RR^n)$
and $f_i$ are not all equal.
Then, in any of cases (a) or (b), there is a metric $\tilde g=(1/\phi^2)g$ solving the problem (\ref{E-PDE-c1}) if and only~if there is a differentiable function $U(x_k)$ such that $\phi=e^{U}$ and (\ref{E-T4-fh-a}) holds.
If~$\phi\le C$ for some constant $C > 0$, then the metrics are complete on~$\RR^n$.
\end{cor}

\begin{ex}\rm
(i) In the case (a) of Theorem~\ref{T-main-2}, assuming that all functions $f_i,f_\alpha$ are constant, we obtain $U=a x_k+b$ and $\phi=e^{a x_k+b}$, where $a,b\in\RR$.

(ii) Consider the function $U=-x_k^{2m}$ for some fixed $k\le p_1$ and $m\in\NN$.
In conditions of Corollary~\ref{C-main-33}, case (a), we obtain
$f_k=-2m(2m-1)p_2x_k^{2m-2}\le0$, $f_i=-4m^2p_2 x_k^{4m-2}\le0$
and $f_\alpha=-2m x_k^{2m-2}[2m-1+2m(p_1-1)x_k^{2m}]\le0$.
Hence $\widetilde\ric_{i}\le0\ (i=1,2)$. By Corollary~\ref{C-main-33}, the metric $\tilde g$ is complete on $\RR^n$.

(iii) Consider the periodic function $U=\sin x_k$ for some fixed $k\le p_1$.
In both cases of Corollary~\ref{C-main-33}, the metric $\tilde g$ is periodic in all variables, and it can be considered as a complete metric on a cylinder or an $n$-dimensional torus. The mixed scalar curvature of $\tilde g$,
$
 \tilde K_{1,2} = -p_2\,e^{2\/\sin x_k}[\,\sin x_k+(p_1-1)\cos x_k],
$
takes positive and negative values. Case (b), $\phi=e^{\sin x_k}$, can be considered as an example of tensors $T$ defined on a flat torus with a pair of complementary distributions ${\mathcal{D}}_i$, that admits a solution to the case (b) of (\ref{E-PDE-c1}).

(iv) From Theorem~\ref{T-main-2}, with $F=x_n$, we obtain results for a half-space $(\RR^n_+, \bar g)$ with the hyperbolic metric $\bar g_{ij}=\delta_{ij}/x_n^2$.
If $U=-x_n^{2m}$, where $m\in\NN$, then $\bar g=\bar g/\phi^2$ is a complete metric on $\RR^n_+$ and the partial Ricci curvatures are negative, the calculations are similar to~(ii).
\end{ex}

\hskip-1mm
In Theorem~\ref{T-main-3},
 $T_{|{\mathcal{D}}_1}=\sum\nolimits_{\,i\le p_1}\hskip-2pt f_i(x_k,x_\delta)\,dx_i^2$ and
 $T_{|{\mathcal{D}}_2}=\sum\nolimits_{\,\alpha>p_1}\hskip-2pt f_\alpha(x_k,x_\delta)\,dx_\alpha^2$ have fixed indices $k\le p_1$ and $\delta>p_1$.

\begin{thm}\label{T-main-3}
Let $(M^n,\bar g)$ be a locally conformally flat Riemannian manifold with complementary orthogonal distributions ${\mathcal{D}}_{1},{\mathcal{D}}_{2}$ of dimensions $p_1,p_2\ge3$.
Let $V\subset M^n$ be an open set with coordinates $(x_1,\ldots,x_n)$ such that
$\bar g_{ij}=\delta_{i j}/F^2$. Given $k\le p_1$ and $\delta>p_1$, consider a symmetric $(0,2)$-tensor $T$ with the properties
 $
 T_{ij}=f_i(x_k,x_\delta)\delta_{ij},\
 T_{\alpha\beta}=f_\alpha(x_k,x_\delta)\delta_{\alpha\beta}$,
 and
 $T_{i\alpha}=0\ (i,j\le p_1,\ \alpha,\beta>p_1),
 $
where the functions $f_i,f_\alpha\in C^1(V)$,
moreover, $f_i$ are not all equal and $f_\alpha$ are not all equal.
Then, in any of cases (a) or (b), there is a metric $\tilde g=\bar g/\phi^2$
solving the problem (\ref{E-PDE-c1}) if and only if there are differentiable functions $v(x_k),w(x_\delta)$ such that $\phi F=v+w$, where
 \begin{eqnarray}\label{E-Rijabc46}
\nonumber
 a) &&\hskip-5mm f_k = [(p_2{v}''+{w}'')(v+w)-p_2\,({v'}^2+{w'}^2)]/(v+w)^2,\\
\nonumber
 &&\hskip-5mm f_\delta = [({v}''+p_1{w}'')(v+w)-p_1\,({v'}^2+{w'}^2)]/(v+w)^2,\\
 &&\hskip-5mm f_\alpha=f_\delta{-}p_1 {w}''/(v+w)\ (\forall\,\alpha\ne\delta),\
    f_i=f_k{-}p_2 {v}''/(v+w)\ (\forall\,i\ne k),\qquad
\end{eqnarray}
 \begin{eqnarray}
\label{E-Rijabc46b}
\nonumber
 b) &&\hskip-5mm f_k =[\,\frac12((2{-}p_2) {v}''+p_1 {w}'')(v+w)
      +p_2(p_1{-}1)\,({v'}^2+{w'}^2)]/(v+w)^2,\\
\nonumber
    &&\hskip-5mm f_\delta =[\,\frac12(p_2 {v}''-(p_1{-}2) {w}'')(v+w)
      +p_1(p_2{-}1)\,({v'}^2+{w'}^2)]/(v+w)^2,\\
    &&\hskip-5mm f_\alpha=f_\delta-p_1 {w}''/(v{+}w)\ (\forall\,\alpha\ne\delta),\
       f_i=f_k-p_2 {v}''/(v{+}w)\ (\forall\,i\ne k).\qquad
\end{eqnarray}
\end{thm}

\begin{thm}\label{T-main-4}
Let $(M^n,\bar g)$ be a locally conformally flat Riemannian manifold with complementary orthogonal distributions ${\mathcal{D}}_{1},{\mathcal{D}}_{2}$ of dimensions $p_1,p_2\ge3$.
Let $V\subset M^n$ be an open set with coordinates $(x_1,\ldots,x_n)$ such that
$\bar g_{ij}=\delta_{i j}/F^2$.
Consider a non-diagonal symmetric $(0,2)$-tensor $T$ with the properties
 $
 T_{ij}=f_{ij}(x_i,x_j),\,
 T_{\alpha\beta}=f_{\alpha\beta}(x_\alpha,x_\beta)$,
 and
 $T_{i\alpha}=0
 \ (i,j\le p_1,\ \alpha,\beta>p_1),
 $
where $f_{AB}\in C^1(V)$.
Suppose that the functions $f_{ij},f_{\alpha\beta}\in C^1(V)$, moreover,
$f_i$ are not all equal and $f_\alpha$ are not all equal.
Then, there is a metric $\tilde g=\bar g/\phi^2$ solving the problem (\ref{E-PDE-c1})(a) if and only~if
up to a change of order of ${\mathcal{D}}_1$ and ${\mathcal{D}}_2$, one of the following cases occur:

{\rm (i)} $f_{12}(x_1,x_2)$ is any nonzero differentiable function, $f_{ij}\equiv0$ for all $i\ne j$
such that $i\ge3$ or $j\ge3$ and $\phi\,F=\varphi(x_1,x_2)$ is a non-vanishing solution to the PDE $\varphi_{,x_1x_2}=(f_{12}/p_2)\varphi$.

{\rm (ii)} There is an integer $p\in[3,p_1]$ such that $f_{ij}=0$, if $i\ne j$, $i\ge p+1$ or $j\ge p+1$.
 Moreover, there exist non-constant differentiable functions, $U_j(x_j)$, for $1\le j\le p$
 such that for all $i,j$, $1\le i\ne j\le p$
 one of the following holds:
\begin{eqnarray}\label{E-P06}
  f_{ij}\eq  p_2 U_i'\, U_j'\quad{\rm and}\quad
  \phi\,F= a\,e^{\sum_{j=1}^p U_j} + b\,e^{-\sum_{j=1}^p U_j},\\
 \label{E-P07}
  f_{ij}\eq -p_2 U_i'\, U_j'\quad{\rm and}\quad
  \phi\,F= a \cos(\sum\nolimits_{j=1}^p U_j) + b \sin(\sum\nolimits_{j=1}^p U_j),
\end{eqnarray}
where $a,b\in\RR$ and $a^2+b^2>0$. Moreover, in each case $\phi$ is defined on an open connected subset of $\,V$, where it does not vanish.
(The solution to (\ref{E-PDE-c1})(b) is constructed similarly.)
\end{thm}

\begin{rem}\label{R-main4-1}\rm
 If $(M^n,\bar g)$ is the Euclidean space, and $|v(x_k)|,\,|w(x_\delta)|\le C$ and $0<|F\phi(x)|\le C$ for some constant $C > 0$, then the metrics given in Theorems~\ref{T-main-3} and \ref{T-main-4} are complete on $\RR^n$.
\end{rem}

By considering $u=(\phi F)^{\,1/(p-1)}$ (when $p_1=p_2=p$) and the mixed scalar curvature $\tilde K$ obtained from the partial Ricci tensor $T_{|\,{\mathcal{D}}_i}$), as a consequence of Theorems~\ref{T-main-1}\,--\,\ref{T-main-4}, case (a), we present $C^\infty$ solutions to the non-linear PDE's of the type (\ref{E-PDE-RnKK}).
We~will show that for certain functions $\tilde K$, depending on functions of one variable, or an arbitrary constant, there exist conformally flat metrics $\tilde g$, whose mixed scalar curvature is $\tilde K$,
see ($\bm P_2$).

\begin{cor}\label{C-main1}
Let $p_1=p_2$ and $\widetilde K$ is defined by

\noindent\ \ \
{\rm (i)} $-p_1p_2[\lambda+2(a_2-a_1)\,\mu]$, where $a_1, a_2, b_k, c\in\RR$, and $\lambda,\,\mu$ are given in~(\ref{E-lambda-mu}).

\noindent\ \
{\rm (ii)} $p_2\,e^{2\/U}[U''-(p_1{-}1){U'}^2]$, where $U(x_k)$ is a differentiable function,
 for some $k\le p_1$.

\noindent\ \
{\rm (iii)} $(v+w)(p_2 {v}'' + p_1 {w}'') - p_1 p_2 ({v'}^2+{w'}^2)$, where
$v(x_k),\,w(x_\delta)$ are differentiable functions, for some $k\le p_1,\,\delta> p_1$.

\noindent\ \
{\rm (iv)} $p_2(a f- b f^{-1})\big[\sum\nolimits_j( {U_j'}^2+U_j'')
 +p_1\frac{a f- b f^{-1}}{a f+ b f^{-1}}\sum\nolimits_j{U_j'}^2\big]$,
 where $U_j(x_j)\ (1\le j\le p)$ are non-constant differentiable functions,
 $3\le p\le p_1,\ a^2+b^2>0$ and $f=e^{\sum_j U_j}$.

Then (\ref{E-PDE-RnKK}) has a solution, globally defined on $\RR^n$, given by

\noindent\ \
{\rm (i)} $u=\big[\sum\nolimits_{i\le p_1}(a_1 x_i^2+b_i x_i)
 +\sum\nolimits_{\alpha>p_1}(a_2 x_\alpha^2+b_\alpha x_\alpha)+c\,\big]^{2/(n-2)}$.

\noindent\ \
{\rm (ii)} $u=e^{\,2 U/(n-2)}$.

\noindent\ \
{\rm (iii)} $u=(v+w)^{2/(n-2)}$.

\noindent\ \
{\rm (iv)} $u=(a f+ b f^{-1})^{2/(n-2)}$.
(If $a=1$ and $b=0$, then
 $\widetilde K= p_2 e^{\sum_j U_j}\big[(p_1+1)\sum\nolimits_j {U_j'}^2+\sum\nolimits_j U_j''\big]$ and (\ref{E-PDE-RnKK}) has a solution $\,u = e^{\,2\,(\sum_j U_j)/(n-2)}$.
 A~solution to (\ref{E-PDE-RnKK}) corresponding to (\ref{E-P07}) is constructed similarly).
\end{cor}

\section{The variational formulae for the total mixed scalar curvature}
\label{sec:scalar}

Now let ${\mathcal D}_1$ and ${\mathcal{D}}_2={\mathcal D}_1^\perp$ be a pair
of complementary orthogonal distributions on a closed Riemannian manifold $(M, g)$.

\begin{df}\label{D-quasip}\rm
Let $p_1\le p_2$, and $\delta_{ij}$ the Kronecker symbol. For any point $x\in M$ and orthonormal bases $e_i\ (i\le p_1)$ of ${\mathcal{D}}_1(x)$ and $\eps_j\ (j\le p_2)$ of ${\mathcal{D}}_2(x)$, consider the bilinear form
${\cal I}_{q}(\vec\omega,\vec\omega)=\sum_{i,j=1}^{p_1}\Phi_{ij}\,\omega_i\,\omega_j$ with the coefficients
\begin{equation}\label{E-phi-ij}
\begin{array}{c}
 \Phi_{ij}=[(\ric_{1}-\ric_{2})(\eps_j,\eps_j)-(\ric_{1}-\ric_{2})(e_i,e_i)]\,\delta_{ij}\\
 +2[\,g(R(e_i,e_j)\eps_i, \eps_j)+g(R(e_i,\eps_j)\eps_i, e_j)].
\end{array}
\end{equation}
We say that ${\cal I}$ is \textit{quasi-positive} if ${\cal I}_{x}$ is positive definite
for arbitrary adapted orthonormal basis $\{e_i,\eps_j\}$ at any $x\in M\setminus\Sigma$,
where $\Sigma$ is a set of zero volume.
\end{df}

In next theorem, using the partial Ricci curvature, we calculate the first and second variations of the \textit{total mixed scalar curvature}
$I_K: {\mathcal D}_1\to\int_M K_{{\mathcal D}_1,{\mathcal D}_2}\,{\rm d}\vol$
of a distribution of arbitrary dimension $p_1,\ 1\le p_1<\dim M$.

\begin{thm}\label{T-Kmix}
 A distribution ${\mathcal D}$ on a closed Riemannian manifold $(M,g)$ is a critical point for the functional $I_K$ if and only if
\begin{equation}\label{eulerK}
 (\ric_{1}-\ric_{2})({\mathcal D}_1,{\mathcal D}_2)=0.
\end{equation}
It is a point of local minimum if the form ${\cal I}$ (of Definition~\ref{D-quasip}) is quasi-positive.
\end{thm}

Let $\dim\,{\mathcal D}_2=1$, and $N\in{\mathcal D}_2$ be a unit vector field on a domain $V\subset M$.
Then $\ric_1({\mathcal D}_1,{\mathcal D}_2)=0$ and (\ref{eulerK}) is reduced to $\ric(X,N)=0\ (X\in{\mathcal D}_1)$.
Only one term of ${\cal I}$ is presented:
$\Phi_{11}=\ric(\eps_1,\eps_1)-\ric(N,N)$. Hence, ${\cal I}$ is positive definite
at $x\in M$ if and only if $\ric(X,X)-\ric(N,N)|X|^2>0$ for all non-zero $X\perp N$
in $T_xM$. With this remark, we have the following.

\begin{cor}[\cite{rw4}]\label{cor:m2}
 A unit vector field $N$ orthogonal to a codimension-one distribution ${\mathcal D}_1$ on a compact Riemannian manifold $(M, g)$ is a critical point for the functional
 $I_2: N\to\int_M\ric(N,N)\,{\rm d}\vol$, if and only~if
 \begin{equation}\label{euler}
 {\ric}(N,X) = 0\qquad \forall\, X\in{\mathcal D}_1.
\end{equation}
It is a point of local minimum if the form ${\cal I}_{2, N}$ is positive definite on the space of sections of ${\mathcal D}_1$. The above bilinear form on the space of vector fields orthogonal to $N$ is given by
${\cal I}_{2,N}(X,Y) = \ric(X, Y) -\ric(N,N)\,g(X, Y)$.
\end{cor}

\begin{ex}\rm
(a) If ${\mathcal D}_1$ and ${\mathcal D}_2$ are curvature invariant, then $\ric_{i}({\mathcal D}_1,{\mathcal D}_2)=0$ $(i=1,2)$, hence ${\mathcal D}_1$ is critical for~$I_K$. A distribution ${\mathcal D}_1$ is called \textit{curvature invariant} if $R(X,Y)Z\in{\mathcal D}_1$ for all $X,Y,Z\in{\mathcal D}_1$.

(b) For Einstein manifold $M^4$ of non-constant sectional curvature, any ``optimally placed"
two-dimensional distribution consists of planes with maximal or minimal curvature.
 To see this, one may use the following characteristic property of Einstein 4-manifolds among Riemannian manifolds $(M^4,g)$: ''the sectional curvature $K(Q)=K(Q^\perp)$ for any 2-plane $Q$",
 see~\cite[Corollary~1.129]{bes}.
(Hence, the product $S^2(1)\times\HH^2(-1)$ is not Einstein manifold, namely, $\ric(e_1,e_1)=1,\,\ric(e_3,e_3)=-1$ when $e_1\in TS^2,\,e_3\in T\HH^2$.)
\end{ex}

The \textit{co-nullity tensor} $C:{\mathcal D}_2\times{\mathcal D}_1\to{\mathcal D}_1$ of a distribution
 ${\mathcal D}_1$ assigns to a~pair $(N, X)$, $X$ being tangent and $N$ normal to ${\mathcal D}_1$,
 the tangent (to ${\mathcal D}_1$) component of the vector field $\nabla_X N$.
Given $N\in{\mathcal D}_2$, let $\sigma_i(C(N,\cdot))\ (0\le i\le p_1)$ be the coefficients of the polynomial
$\det({\rm Id}+t C(N,\cdot))=\sum_i\sigma_i(C(N,\cdot))\,t^i$.

The \textit{$2k$-th mean curvature\/} of a distribution ${\mathcal D}_1^{p_1}$ is the integral
\begin{equation*}
 \sigma_{2k}({\mathcal D}_1(x))=\frac {p_1}{\vol(S^{p-1})}\int_{\,N\perp{\mathcal D}_1(x),\,|N|=1}
 \hskip-3mm\sigma_{2k}(C(N,\cdot))\,{\rm d}\,\omega, \quad \mbox{for all }\ x\in M.
\end{equation*}
The same formula determines $\sigma_{2k-1}({\mathcal D}(x))=0$.

 Denote by
$C_{1,ij}^{\,\alpha}{=}g(C_{1}^{\,\alpha}e_i,e_j)$
and
$C_{2,\alpha\beta}^{\,i}{=}g(C_{2}^{\,i}e_\alpha,e_\beta)$,
where the linear operators
$C_{1}^{\,\alpha}:\RR^{p_1}\to\RR^{p_1}$
and
$C_{2}^{\,i}:\RR^{p_2}\to\RR^{p_2}$
correspond to co-nullity tensors $C_i$ of ${\mathcal{D}}_i$ at $x\in M$.

\begin{prop}\label{L-IntK}
Let ${\mathcal{D}}_i,\,\dim {\mathcal{D}}_i=p_i\ (i=1,2)$ be a pair of complementary orthogonal distributions on a compact Riemannian manifold $(M,g)$.
Then
$$I_K=2\int_M\big(\sigma_{2}({\mathcal{D}}_1)+\sigma_{2}({\mathcal{D}}_2)\big)\,{\rm d}\vol.$$
\end{prop}

 The extremal values of $I_K$ can be used for estimation
 of the total energy and bending of a distribution or a vector field.

We can regard the distribution ${\mathcal{D}}_i$ as the map $\tilde {\mathcal{D}}_i:M\to G(p_i,M)$
(section of the Grassmann bundle $G(p_i,M)=\cup_{x\in M} G(p_i,T_x M)$),
where $\tilde {\mathcal{D}}_1(x)=e_{1}\wedge\dots\wedge e_{p_1}$ and $\tilde {\mathcal{D}}_2(x)=e_{p_1+1}\wedge\dots\wedge e_{m}$ are the $p_i$-vectors determined locally by ${\mathcal{D}}_1(x)$ and ${\mathcal{D}}_2(x)$, resp.
For a map between Riemannian spaces $f: \bar M\to (M,g)$, the \textit{energy} is defined to be
 $\mathcal{E}(f)=\frac12\int_{M}\sum_{\,a=1}^m g(d\,f(e_a), d\,f(e_a))\,d\vol$, see~\cite{el}.
The \textit{corrected energy} of a $p_2$-dimensional distribution ${\mathcal{D}}_2$ on a $(p_1+p_2)$-dimensional Riemannian manifold $(M, g)$ is defined in \cite{c} as
\begin{equation*}
 \mathcal{D}({\mathcal{D}}_2)=
 \int_{M}\sum\nolimits_{\,a=1}^m
 \Big[|\nabla_{e_a}\tilde {\mathcal{D}}_2|^2 +p_1(p_1-2)\,|H_1|^2+p_2^{\,2}\,|H_2|^2\Big]{\rm d}\vol,
\end{equation*}
where $|d\,\tilde {\mathcal{D}}_2|$ is calculated from the definition of Sasaki metric $g_s$:
\begin{equation*}
 \sum\nolimits_{\,a=1}^m g_s(d\,\tilde{\mathcal{D}}_2(e_a), d\,\tilde{\mathcal{D}}_2(e_a))=
 \sum\nolimits_{\,a=1}^m [\,g(e_a,e_a)+g(\nabla_{e_a}\tilde{\mathcal{D}}_2, \nabla_{e_a}\tilde {\mathcal{D}}_2)]
\end{equation*}
and $H_1=-\frac1{p_1}\sum_{\,\alpha}(\sum_i C_{1,\,i\,i}^{\,\alpha})\,e_\alpha$,
    $H_2=-\frac1{p_2}\sum_{\,i}(\sum_\alpha C_{2,\,\alpha\,\alpha}^{\,i})\,e_i$.
If ${\mathcal{D}}_2$ is integrable, then
 $\mathcal{D}({\mathcal{D}}_2)\ge\int_M K_{1,2}\,{\rm d}\vol$, see \cite{c}.
Similarly, we define the \textit{total bending}
 $\mathcal{B}({\mathcal{D}}_2)=c_{n}\int_M |\nabla\tilde {\mathcal{D}}_2|^2\,{\rm d}\vol$,
where $c_{n}$ is a constant.

\begin{prop}\label{P-bend} The total bending of a $p_1$-dimensional distribution ${\mathcal{D}}_1$ on a
$(p_1+p_2)$-dimensional Riemannian manifold $M$ satisfies the inequality
\begin{equation}\label{E-bendD5}
 \mathcal{B}({\mathcal{D}}_1)\ge c_{n}\int_M
 \Big(\,\frac2{p_1{-}1}\sum\nolimits_{\,\alpha}\sigma_2(C_{1}^{\,\alpha})
 +\frac2{p_2{-}1}\sum\nolimits_{\,i}\sigma_2(C_{2}^{\,i})\Big)\,d\vol
\end{equation}
that for $p_1=p_2=p$ takes the form $\mathcal{B}({\mathcal{D}}_1)\ge\frac{c_{n}}{p-1}\,I_K$.
\end{prop}

\section{Proof of results}
\label{sec:proof}

\noindent\textbf{Proof of Proposition \ref{P-main1}}.
The partial Ricci curvature in local coordinates~is
\begin{equation*}
\begin{array}{c}
 \ric_{1}(g)_{ij}=\frac12\sum_{\alpha\beta}g^{\alpha\beta}
 (g_{i\alpha,j\beta}+g_{j\alpha,i\beta}-g_{ij,\alpha\beta}-g_{\alpha\beta,ij})
 +Q_1(g,\partial g),\\
  \ric_{2}(g)_{\alpha\beta}=\frac12\sum_{ij}g^{ij}
 (g_{\alpha i,\beta j}+g_{\beta i,\alpha j}-g_{\alpha\beta, ij}-g_{ij, \alpha\beta})
 +Q_2(g,\partial g),
\end{array}
\end{equation*}
where $g^{a b}$ is the inverse of $g_{a b}$ and $Q_i$ is a function of $g$ and its derivatives, and is homogeneous of degree 2 in the first derivatives of $g$.
 In view of $g_{i\alpha}\equiv0$, we~obtain
\begin{equation}\label{E-PDE5}
\begin{array}{c}
 \ric_{1}(g)_{ij}=-\frac12\sum_{\alpha\beta}g^{\alpha\beta}
 (g_{ij,\alpha\beta}+g_{\alpha\beta,ij})+Q_1(g,\partial g),\\
  \ric_{2}(g)_{\alpha\beta}=-\frac12\sum_{ij}g^{ij}
 (g_{\alpha\beta, ij}+g_{ij, \alpha\beta})+Q_2(g,\partial g).
\end{array}
\end{equation}
Assume that $g$ has the form
$g = \sum_{\,a=1}^{n}\big(1-\sum_{\,b\le n} c_{bb} x_b^2\big) dx^a\otimes dx^a$.
Then
 $g_{ii,\alpha\alpha}=2\,c_{\alpha\alpha},\,g_{\alpha\alpha,ii}=2\,c_{ii}$.
Substituting in (\ref{E-PDE5}) and using (\ref{E-PDE-c1}), we get at~$O$
 \begin{equation}\label{E-PDE6}
 p_2 c_{i i}{+}\sum\nolimits_\alpha c_{\alpha\alpha}= T_{ii}\ (1\le i\le p_1),\quad
 p_1 c_{\alpha\alpha}{+}\sum\nolimits_i c_{i i}=T_{\alpha\alpha}\ (p_1<\alpha\le n).
\end{equation}
The summing of first $p_1$ equations and last $p_2$ ones in (\ref{E-PDE6}) yields
 $p_1\sum_\alpha c_{\alpha\alpha}+p_2\sum_i c_{i i}=K$,
where
$K=\sum_{i\le p_1} T_{ii}=\sum_{\alpha>p_1} T_{\alpha\alpha}$.
The linear system (\ref{E-PDE6}) consists of $n$ equations and the same number of variables, its rank $<n$.
It~is easily seen that, for instance, $c_{ii}=T_{ii}/p_2\ (1\le i\le p_1)$ and $c_{\alpha\alpha}=(T_{\alpha\alpha}-K/p_2)/p_1\ (p_1<\alpha\le n)$ is a solution (satisfying
$\sum_\alpha c_{\alpha\alpha}=0$).
The~tensor $g$ is positive definite at $O$ in a small neighborhood.
\hfill$\square$

\vskip1mm
\textbf{Proof of Proposition~\ref{E-RicD}}.
Define a metric $\tilde g$ on $M$ by $\tilde g=e^{\psi} g$,
where $\psi=-2\log\phi$. The connections $\nabla$ and $\tilde\nabla$ of $g$ and $\tilde g$ are related~by \cite{GKM}
 $\tilde\nabla_X Y=\nabla_X Y+\frac12(X(\psi)\,Y + Y(\psi)\,X -g(X, Y)\nabla\psi).$
Let $\tilde R$ be the curvature tensor of $\tilde g$.
Using known formulae (see, for example, \cite{GKM}) we obtain
\begin{eqnarray}\label{E-R}
\nonumber
 &&\hskip-5mm\tilde R(X,Y)Y = R(X,Y)Y
 +\frac12(g(h_\psi(X),Y)Y-g(h_\psi(Y),Y)X-|Y|^2 h_\psi(X))\\
 &&\hskip-6mm+\frac14([Y(\psi)^2{-}|Y|^2|\nabla\psi|^2]X
 {+}[X(\psi)Y(\psi){-}Y(\psi)|\nabla\psi|^2]Y{+}X(\psi)|Y|^2\nabla\psi)
\end{eqnarray}
where $X\in {\mathcal{D}}_{2}, Y\in {\mathcal{D}}_{1}$, hence $g(X,Y)=0$.
Note that if $Z$ is a $g$-unit vector then $\tilde Z=Z\,e^{-\psi/2}$ is a $\tilde g$-unit vector.
From this and (\ref{E-R}), using adapted orthonormal base $\{e_i, e_\alpha\}$, one may deduce the relation between partial Ricci curvatures in both metrics
\begin{eqnarray*}
 \widetilde\ric_{1}(e_A,e_B)\eq\ric_{1}(e_A,e_B)-\frac12(\delta_{AB}\sum\nolimits_\alpha g(h_\psi(e_\alpha),e_\alpha)+p_2 g(h_\psi(e_A),e_B))\\
 \plus\frac14([\sum\nolimits_\alpha g(\nabla\psi, e_\alpha)^2 -p_2|\nabla\psi|^2]\delta_{AB}+p_2 g(\nabla\psi, e_A)\,g(\nabla\psi, e_B)).
\end{eqnarray*}
In matrix notation, this reads as
 \begin{equation}\label{E-Rij2}
 \widetilde\ric_{1} = \ric_{1} -\frac12\big((\Delta^{(2)}\psi)\,g +p_2 h_\psi\big)
 {+}\frac14\big(\big[|\nabla^{(2)}\psi|^2-p_2|\nabla\psi|^2\big]g
  +p_2 g(\nabla_{\cdot}\psi, \nabla_{\cdot}\psi)\big)
\end{equation}
where $\nabla^{(2)}$ and $\Delta^{(2)}$ are ${\mathcal{D}}_{2}$-gradient and ${\mathcal{D}}_{2}$-laplacian of a function.

To shorten the formulae, we will turn back to the function $e^\psi=1/\phi^2$.
In this case
 $\nabla\psi=-\frac2{\phi}\nabla\phi$,
 $h_\psi(e_A,e_B)=\frac2{\phi^2}\,g(\nabla_{e_A}\phi, \nabla_{e_B}\phi)
 -\frac2{\phi}\,h_\phi(e_A,e_B)$,
 $\Delta^{(2)}\psi=\frac2{\phi^2}\,|\nabla^{(2)}\phi|^2-\frac2{\phi}\,\Delta^{(2)}\phi$,
etc. Substituting above equalities in (\ref{E-Rij2}), we obtain
\begin{eqnarray*}
 \widetilde\ric_{1}-\ric_{1}\eq
 -\big(\frac1{\phi\,^2}\,|\nabla^{(2)}\phi|^2 -\frac1{\phi}\,\Delta^{(2)}\phi\big)g
 -p_2\big(\frac1{\phi\,^2}\,g(\nabla_{\cdot}\phi,\nabla_{\cdot}\phi) -\frac1{\phi}\,h_\phi\big)\\
 \plus\frac1{\phi\,^2}|\nabla^{(2)}\phi|^2 g -\frac{p_2}{\phi\,^2}\,|\nabla\phi|^2 g
 +\frac{p_2}{\phi\,^2}\,g(\nabla_{\cdot}\phi, \nabla_{\cdot}\phi)
\end{eqnarray*}
that is simplified to (\ref{E-Rij3a}). The formula for ${\mathcal{D}}_2$ is proved similarly.
By
$
 \widetilde K_{1,2}=\phi\/^2\sum\nolimits_{\,i\le p_1}\hskip-3pt\widetilde{\ric}_{1}(e_i,e_i),
$
(\ref{E-tildeK1}) is the result of the trace operation applied to (\ref{E-Rij3a}).\hfill$\square$

\vskip1mm
\textbf{Proof of Corollary~\ref{Cor-T1-u}}.
Recall that $\tilde g=(1/\phi^2)\,g$.
Let $\phi=u^{-\gamma}$.
 Using $\phi=u^{-\gamma}$
and $\Delta^{(i)}={\rm Div}^{(i)}\nabla^{(i)}$, we find
 $|\nabla\phi|^2=\gamma^2 u^{-2\gamma-2}|\nabla u|^2$
and
$
 \Delta^{(i)}\phi=-\gamma u^{-\gamma-1}\Delta^{(i)}u
 +\gamma(\gamma+1)\,u^{-\gamma-2}|\nabla^{(i)} u|^2\ (i=1,2)$.
Substituting  in (\ref{E-tildeK1}), we obtain
\begin{eqnarray*}\label{E-PDE-RnKK-first}
 \bar K u^{2\gamma+1}\eq K u-\gamma\big(p_1\Delta^{(2)}+ p_2\Delta^{(1)}\big)u\\
 \plus\!\gamma\,u^{-1}\big[(\gamma+1)(p_1|\nabla^{(2)} u|^2\!+ p_2|\nabla^{(1)} u|^2\big)
 \!-\gamma\,p_1p_2|\nabla u|^2\big].
\end{eqnarray*}
For $n$ even, $p_1=p_2=n/2$ and $\gamma=2/(n-2)$ this yields~(\ref{E-PDE-RnKK}).\hfill$\square$

\vskip1mm
\textbf{Proof of Theorem~\ref{T-main-1}}.
a) The compatibility condition for ($\bm P_1$) is $p_1 f_1=p_2 f_2$.
Observe that $\tilde g=\bar g/\phi^2=g/(\phi F)^2=g/{\varphi}^2$, where $g$ is the Euclidean metric, and ${\varphi}=\phi F$. In view of $\ric_{1}=\ric_{2}=0$ for $g$, we~have, see~(\ref{E-Rij3a}),
 \begin{eqnarray}
 \label{E-Rij3aH}
 \nonumber
 \widetilde\ric_{1} = \big[p_2{\varphi}\,h_{{\varphi}}
 +\big({\varphi}\,\Delta^{(2)}{\varphi}-p_2\,|\nabla{\varphi}\,|^2\big)g\big]/{{\varphi}\,^2},\\
 \widetilde\ric_{2} = \big[p_1{\varphi}\,h_{{\varphi}}
 +\big({\varphi}\,\Delta^{(1)}{\varphi}-p_1\,|\nabla{\varphi}\,|^2\big)g\big]/{{\varphi}\,^2}.
\end{eqnarray}
Since $T_{|{\mathcal{D}}_k}=\widetilde\ric_{k\,| {\mathcal{D}}_k}\ (k=1,2)$, we obtain
\begin{eqnarray*}
 {{\varphi}\,^2}f_1 g = p_2{\varphi}\,h_{{\varphi}}
 +\big({\varphi}\,\Delta^{(2)}{\varphi}-p_2\,|\nabla{\varphi}|^2\big)g
 \quad{\rm on}\ {\mathcal{D}}_{1},\\
 {{\varphi}\,^2}f_2 g = p_1{\varphi}\,h_{{\varphi}}
 +\big({\varphi}\,\Delta^{(1)}{\varphi}-p_1\,|\nabla{\varphi}|^2\big)g
 \quad{\rm on}\ {\mathcal{D}}_{2}.
\end{eqnarray*}
Hence, the problem is reduced to studying the following system of PDE's:
\begin{eqnarray}\label{E-Rij4abcH}
\nonumber
 p_2{\varphi}_{,x_i x_i}\eq{\varphi} f_1 -\Delta^{(2)}{\varphi}+{p_2}\,|\nabla{\varphi}\,|^2/{\varphi},\quad
 i\le p_1\\
 p_1{\varphi}_{,x_\alpha x_\alpha}\eq{\varphi} f_2 -\Delta^{(1)}{\varphi}+{p_1}\,|\nabla{\varphi}\,|^2/{\varphi}, \quad \alpha>p_1,\\
\nonumber
 {\varphi}_{,x_k x_m}\eq0,\quad 1\le k\ne m\le n.
\end{eqnarray}
From the last equation of (\ref{E-Rij4abcH}) we conclude that
${\varphi}=\sum_{k=1}^{n}\phi_k(x_k)$.
From the first two equations of (\ref{E-Rij4abcH}) we deduce $\phi_i''(x_i)=2 a_1\in\RR$ for all $i\le p_1$ and $\phi_\alpha''(x_\alpha)=2 a_2\in\RR$ for all $\alpha>p_1$.
Therefore,
 ${\varphi}=\sum_{\,i\le p_1}(a_1 x_i^2+b_i x_i)+
 \sum_{\,\alpha>p_1}(a_2 x_\alpha^2+b_\alpha x_\alpha)+c$,
 where $b_k, c\in\RR$. We also have $\Delta^{(1)}{\varphi}=2a_1 p_1$ and $\Delta^{(2)}{\varphi}=2a_2 p_2$.
 Hence the first two equations of (\ref{E-Rij4abcH}) are reduced~to
 \begin{equation*}
 2p_2(a_1+a_2)={\varphi} f_1 +{p_2}\,|\nabla{\varphi}|^2/{\varphi},\quad
 2p_1(a_1+a_2)={\varphi} f_2 +{p_1}\,|\nabla{\varphi}|^2/{\varphi}.
\end{equation*}
Comparing them, we see that the equality $p_1f_1=p_2f_2\,(=\phi^2\widetilde K_{1,2})$ is necessary for the solution existence. In view of $|\nabla{\varphi}|^2-2(a_1+a_2){\varphi}=\lambda-2(a_2-a_1)\,\mu$, we obtain $f_1$ and $f_2$, as required.

If $a_1a_2\le0$ and $a_1^2+a_2^2+\sum_k b_k^2>0$ then the set $\{\phi=0\}$
of singularities of $\tilde g$ is non-empty and can be explicitly described.
If $a_1a_2>0$ then the inequality $\phi>0$ means that the discriminant of a quadratic equation is negative.

b) The proof is similar to the previous one.
The compatibility condition for ($\bm P_2$) is $(1-p_2/2)p_1 f_1=(1-p_1/2)p_2 f_2$.
The problem is reduced to the
PDE's:
\begin{eqnarray}\label{E-Rij5abcH}
\nonumber
 p_2{\varphi}_{,x_i x_i}\eq{\varphi} f_1 -\Delta^{(2)}{\varphi}+{p_2}\,|\nabla{\varphi}|^2/{\varphi}
 +\widetilde K_{1,2}/(2{\varphi}),\quad i\le p_1\\
 p_1{\varphi}_{,x_\alpha x_\alpha}\eq{\varphi} f_2 -\Delta^{(1)}{\varphi}+{p_1}\,|\nabla{\varphi}|^2/{\varphi}
 +\widetilde K_{1,2}/(2{\varphi}),\quad \alpha>p_1,\\
\nonumber
 {\varphi}_{,x_k x_m}\eq0,\quad 1\le k\ne m\le n,
\end{eqnarray}
where ${\varphi}=F\phi$.
From the last equation of (\ref{E-Rij5abcH}) we conclude that ${\varphi}=\sum_{k=1}^{n}\phi_k(x_k)$.
Moreover, from the first two equations of (\ref{E-Rij5abcH}) we deduce
${\varphi}_i''(x_i)=2 a_1\in\RR$ for all $i\le p_1$
and
${\varphi}_\alpha''(x_\alpha)=2 a_2\in\RR$ for all $\alpha>p_1$.
Therefore, ${\varphi}=\sum_{\,i\le p_1}(a_1 x_i^2+b_i x_i)+
 \sum_{\,\alpha>p_1}(a_2 x_\alpha^2+b_\alpha x_\alpha)+c$, where $b_k, c\in\RR$.
 We also calculate
$\Delta^{(1)}{\varphi}=2a_1 p_1$ and $\Delta^{(2)}{\varphi}=2a_2 p_2$.
 Hence the first two equations of (\ref{E-Rij5abcH}) are reduced to the equations
\begin{eqnarray*}
 2p_2(a_1+a_2){\varphi}\eq{\varphi}^2 f_1 +{p_2}\,|\nabla{\varphi}|^2
 +\widetilde K_{1,2}/2,\\
 2p_1(a_1+a_2){\varphi}\eq{\varphi}^2 f_2 +{p_1}\,|\nabla{\varphi}|^2 +\widetilde K_{1,2}/2.
\end{eqnarray*}
Since $|\nabla{\varphi}|^2-2(a_1+a_2){\varphi}=\lambda-2(a_2-a_1)\,\mu$, we obtain $f_1$ and $f_2$, as required.$\,\square$

\vskip1mm
\textbf{Proof of Corollary~\ref{Cor-main1}}.
a) We set $F=1$ and obtain $\phi, f_1, f_2$ as in the proof of Theorem~\ref{T-main-1}.
 Assume that $a_1=a_2=a$, see Example~\ref{R-a1-a2}. Then we have
 $\phi=\sum\nolimits_{i=k}^{n}(a x_k^2+b_k x_k)+c$, and
 $f_1=-p_2\lambda/\phi^2,\
  f_2=-p_1\lambda/\phi^2,\
 \widetilde K_{1,2}=-p_1p_2\lambda\,$.
 If $\lambda<0$ then $a\ne0$ and $\phi>0$. Namely,
 $\phi\ge c-\frac{1}{4a}\sum\nolimits_k b_k^2 >0$, if $a>0$, and
 $\phi\le c+\frac{1}{4a}\sum\nolimits_k b_k^2 <0$, if $a<0$.
 Let $\lambda=0$. If $a\ne0$ then $\phi=0$ has a unique solution at
 $\tilde x=-(b_1,\ldots, b_n)/(2a)$, hence $\tilde g$ has one singular point.
 If $a=0$ then $b_k=0$. Hence $\phi=c$ and $\tilde g$ reduces to a homothety.

 Now let $\lambda>0$.
If $a=0$ then there is $k_0$ such that $b_{k_0}\ne0$ and $\phi$ vanishes on the hyperplane $(\sum_k b_k x_k) + c = 0$ -- a singularity set of $\tilde g$. If $a\ne0$ then $\phi$ vanishes on the $(n-1)$-dimensional sphere, centered at $\tilde x=-(b_1,\ldots, b_n)/(2a)$ with radius $\sqrt{\lambda}/(2|a|)$ -- a singularity set of $\tilde g$.
b) The proof is similar to the previous one.\hfill$\square$

\vskip1mm
\textbf{Proof of Theorem~\ref{T-main-2}}.
a) We set $\tilde g=\bar g/\phi^2=g/(\phi F)^2=g/{\varphi}^2$, where $g$ is the Euclidean metric, and ${\varphi}=\phi F$. From (\ref{E-PDE-c1})(a), (\ref{E-Rij3a}), in view of $\ric_{1}=\ric_{2}=0$ for $g$, we~have (\ref{E-Rij3aH}). Hence, the problem reduces to studying the following system:
\begin{eqnarray}\label{E-Rijabc4}
\nonumber
 p_2{\varphi}_{,x_i x_i}\eq{\varphi} f_i -\Delta^{(2)}{\varphi}+{p_2}\,|\nabla{\varphi}|^2/{\varphi},
 \quad i\le p_1,\\
 p_1{\varphi}_{,x_\alpha x_\alpha}\eq{\varphi}\,f_\alpha -\Delta^{(1)}{\varphi}+{p_1}\,|\nabla{\varphi}|^2/{\varphi},
  \quad \alpha>p_1,\\
\nonumber
 {\varphi}_{,x_k x_m}\eq0,\quad 1\le k\ne m\le n.
\end{eqnarray}
From the third equation of (\ref{E-Rijabc4}) we conclude that
${\varphi}=\sum_{s=1}^{n}\phi_s(x_s)$, which substituted in the first two equations gives
\begin{eqnarray}\label{E-Rijabc44}
\nonumber
 p_2\phi_i''(x_i)\eq{\varphi} f_i-\Delta^{(2)}{\varphi}+{p_2}\,|\nabla{\varphi}|^2/{\varphi},
 \quad i\le p_1,\\
 p_1\phi_\alpha''(x_\alpha)\eq
 {\varphi}\,f_\alpha-\Delta^{(1)}{\varphi}+{p_1}\,|\nabla{\varphi}|^2/{\varphi},
  \quad \alpha>p_1,
\end{eqnarray}
where $f_i=f_i(x_k),\,f_\alpha=f_\alpha(x_k)$.
As a consequence of (\ref{E-Rijabc44}) we have
\begin{eqnarray}\label{E-Rijabc45-A}
\nonumber
 p_2(\phi''_k(x_k)-\phi''_i(x_i))\eq {\varphi}(f_k-f_i),\quad i\le p_1,\\
 p_1(\phi''_\alpha(x_\alpha)-\phi''_\beta(x_\beta))\eq
  {\varphi}(f_\alpha-f_\beta),\quad \alpha,\beta>p_1,
\end{eqnarray}
and then
\begin{equation}\label{E-Rijabc45}
 \frac{\phi''_k(x_k){-}\phi''_i(x_i)}{\phi''_k(x_k){-}\phi''_j(x_j)}{=}
 \frac{f_k{-}f_i}{f_k{-}f_j}\ (i,j\le p_1),\
 \frac{\phi''_\alpha(x_\alpha){-}\phi''_\beta(x_\beta)}
 {\phi''_\alpha(x_\alpha){-}\phi''_\gamma(x_\gamma)}{=}
 \frac{f_\alpha{-}f_\beta}{f_\alpha{-}f_\gamma} \ (\alpha,\beta,\gamma>p_1).
\end{equation}
From (\ref{E-Rijabc45})$_1$ (not all $f_i$ are equal), in view of $p_1,p_2\ge3$, we deduce
that ${\varphi}$ is the function of $x_k$ only, i.e., ${\varphi}=\phi_k(x_k)$.
Next, from (\ref{E-Rijabc45-A})$_2$ we obtain
 $f_\alpha=f_\beta\ (\forall\,\alpha,\beta)$ and
 $f_i=f_k-p_2\phi''_k/\phi_k\ (\forall\,i\ne k)$.
 Then we calculate $|\nabla{\varphi}|^2={\phi_k'}^2$,
 $\Delta^{(1)}{\varphi}=\phi''_k$ and $\Delta^{(2)}{\varphi}=0$.
 Hence, $f_k=p_2\frac{\phi''_k\phi_k-{\phi_k'}^2}{(\phi_k)^2}=p_2(\log\phi_k)''$.
 Assuming $f_k=p_2 U''(x_k)$ we immediately obtain $\phi_k=e^{U(x_k)}$ and
 $\phi={\varphi}/F=\frac{1}{F}e^{U(x_k)}$ that is required. Next, we confirm that $f_i=-p_2{U'}^2\ (i\ne k)$ and $f_\alpha=U''-(p_1-1){U'}^2$.
 One may verify that the compatibility condition for $T$, $f_k+(p_1-1)f_i=p_2 f_\alpha$ holds,
 is satisfied. Finally, $\widetilde K_{1,2}={\varphi}^2 p_2 f_\alpha =p_2e^{2\/U}[U''-(p_1{-}1){U'}^2]$.

 b) The proof is similar to the previous one. The problem is reduced to studying the system
\begin{eqnarray}\label{E-Rijabc4-5}
\nonumber
 p_2{\varphi}_{,x_i x_i}\eq\phi f_i -\Delta^{(2)}{\varphi}+{p_2}\,|\nabla{\varphi}|^2/{\varphi}
 +\widetilde K_{1,2}/(2{\varphi}),
 \quad i\le p_1,\\
 p_1\phi_{,x_\alpha x_\alpha}\eq{\varphi} f_\alpha-\Delta^{(1)}{\varphi} +{p_1}\,|\nabla{\varphi}|^2/{\varphi}+\widetilde K_{1,2}/(2{\varphi}),
 \quad\alpha>p_1,\\
\nonumber
 \phi_{,x_k x_m}\eq0,\quad 1\le k\ne m\le n.
\end{eqnarray}
 As in the case a) we obtain $\phi=\sum_{s=1}^{n}\phi_s(x_s)$, and
\begin{eqnarray}\label{E-Rijabc55}
\nonumber
 p_2\phi''_i(x_i)\eq{\varphi} f_i(x_k) -\Delta^{(2)}{\varphi}+{p_2}\,|\nabla{\varphi}|^2/{\varphi}
 +\widetilde K_{1,2}/(2{\varphi}),
 \quad i\le p_1,\\
 p_1\phi''_\alpha(x_\alpha)\eq{\varphi} f_\alpha(x_k) -\Delta^{(1)}{\varphi}+{p_1}\,|\nabla{\varphi}|^2/{\varphi}
 +\widetilde K_{1,2}/(2{\varphi}),
  \quad \alpha>p_1.
\end{eqnarray}
As a consequence of (\ref{E-Rijabc55}) we again have (\ref{E-Rijabc45-A}), (\ref{E-Rijabc45}). Similarly to case a), we conclude that ${\varphi}=\phi_k(x_k)$.
Next, from (\ref{E-Rijabc45-A}) we obtain
 $f_\alpha=f_\beta\ (\forall\,\alpha,\beta)$ and
 $f_i=f_k-p_2\phi''_k/\phi_k\ (\forall\,i\ne k)$.
 Then we calculate $|\nabla{\varphi}|^2={\phi_k'}^2$,
 $\Delta^{(1)}{\varphi}=\phi''_k$ and $\Delta^{(2)}{\varphi}=0$.
 Hence,
 $\widetilde K_{1,2}(x)=p_2[\phi''_k\phi_k-p_1{\phi_k'}^2]$.
 From (\ref{E-Rijabc55}) with $i=k$ we obtain
 $$f_k=\frac{p_2}{\phi\,^2_k}\big[({\phi''_k\phi_k-{\phi_k'}^2})
 -({\phi''_k\phi_k-p_1{\phi_k'}^2})/2\big]
 =\frac{p_2}2[(\log\phi_k)''+(p_1-1)(\phi_k'/\phi_k)^2].
 $$
 Assuming $\phi_k=e^{U(x_k)}$ we immediately obtain $f_k=\frac12p_2(U''+(p_1-1){U'}^2)$ that is required.
 Next, we obtain $f_i=-\frac12p_2(U''-(p_1{-}3){U'}^2)$ $(i\ne k)$
 and $f_\alpha=\frac12(p_2-2)[-U''+(p_1-1){U'}^2]$.
 One may verify that the compatibility condition for $T$,
 $(1-p_2/2)[f_k+(p_1-1)f_i]=(1-p_1/2)p_2 f_\alpha$, is satisfied.\hfill$\square$

\vskip1mm
\textbf{Proof of Corollary~\ref{C-main-33}} We consider the function $F=1$ and apply the arguments similar to those of Theorem~\ref{T-main-2}. The metric $\tilde g$,
 satisfying $\phi(x_k)\le C$, is complete, since there is a constant $m > 0$ 
 such that $|v|_{\tilde g}\ge m|v|$ for any 
 $v\in\RR^n.\,\square$

\vskip1mm
\textbf{Proof of Theorem~\ref{T-main-3}}.
a) We set ${\varphi}=F\phi$, and as in the proof of Theorem~\ref{T-main-2}, obtain ${\varphi}=\sum_{s=1}^{n}\phi_s(x_s)$.
Similarly to the proof of Theorem~\ref{T-main-2}, we deduce
(\ref{E-Rijabc44}) -- (\ref{E-Rijabc45}),
where $f_i=f_i(x_k,x_\delta),\,f_\alpha=f_\alpha(x_k,x_\delta)$.
 From these (not all $f_i$ are equal and not all $f_\alpha$ are equal), in view of $p_1,p_2\ge3$,
 we have that ${\varphi}$ is the function of $x_k,x_\delta$ only, i.e., ${\varphi}=v+w$,
where $v=\phi_k(x_k)$ and $w=\phi_\delta(x_\delta)$.
Hence
\begin{equation}\label{E-comp-hf}
 f_\alpha=f_\delta-p_1 {w}''/(v+w)\ \ (\forall\,\alpha\ne\delta),\quad
 f_i=f_k-p_2 v''/(v+w)\ \ (\forall\,i\ne k).
\end{equation}
 Then we find $|\nabla{\varphi}|^2={v'}^2+{w'}^2$, $\Delta^{(1)}{\varphi}=v''$ and $\Delta^{(2)}{\varphi}={w}''$.
 Hence, the functions $v,\,w$ satisfy the system (\ref{E-Rijabc46}).
 The compatibility condition for $T$,
\begin{equation}\label{E-comp4}
 p_2(p_1-1) v''-p_1(p_2-1) {w}''=(p_1f_k-p_2f_\delta)(v+w)
\end{equation}
that is the linear combination of  equations in (\ref{E-Rijabc46}) with coefficients $p_1$ and~$p_2$.

b) As in proof of the case a), we conclude that ${\varphi}=\sum_{s=1}^{n}\phi_s(x_s)$. Similarly to the proof of case a), we deduce (\ref{E-Rijabc44}) -- (\ref{E-Rijabc45}).
As a consequence of these we have, where $f_i=f_i(x_k,x_\delta),\,f_\alpha=f_\alpha(x_k,x_\delta)$.
As in a), we deduce that ${\varphi}$ is the function of $x_k,x_\delta$ only, i.e., ${\varphi}=v+w$,
where $v=\phi_k(x_k),\,w=\phi_\delta(x_\delta)$.
Next, we obtain (\ref{E-comp-hf}).
 Then we calculate $|\nabla{\varphi}|^2={v'}^2+{w'}^2$, $\Delta^{(1)}{\varphi}=v''$ and $\Delta^{(2)}{\varphi}={w}''$. The mixed scalar curvature is
 $\tilde K_{1,2}=(v+w)(p_2 v''+p_1 {w}'')-p_1p_2({v'}^2+{w'}^2)$.
 Hence, the functions $v,\,w$ satisfy the nonlinear system (\ref{E-Rijabc46b}).
 The compatibility condition for $T$ takes the form
 $p_2(p_1-1)(p_1-2) v''-p_1 p_2(p_2-2) {w}''=(p_1(p_1-2)f_k-p_2(p_2-2)f_\delta)(v+w)$
 that is the linear combination of equations in (\ref{E-Rijabc46b}).\hfill$\square$

\vskip1mm
\textbf{Lemma A} \cite{pt09}
{\it Assume $\varphi(x_1,\ldots, x_p),\ p\ge 3$, is a non-vanishing differentiable
function that satisfies a system of equations
\begin{equation}\label{E-P20}
  \varphi_{,ij}=f_{ij}\,\varphi,\quad i\ne j,
\end{equation}
where $f_{ij}=f_{ji}$ is a differentiable function of $x_i$ and $x_j$. Assume there is an
open subset $V\subset\RR^p$, where all $f_{ij}$ do not vanish. Then there is an open dense
subset of $V$ where $\prod_i\varphi_{,i}$ does not vanish. On each connected component of
this subset, there exist differentiable functions $U_i(x_i)\ne0\ (i = 1,\ldots, p)$ such that}
  $f_{ij}=U_i(x_i)\,U_j(x_j),\ 1\le i\ne j\le p$.

\vskip1mm\hskip-3.5mm
\textbf{Lemma B} \cite{pt09}
\textit{A non-vanishing differentiable function $\varphi(x_1,\ldots,x_p),\ p\ge 3$,
\begin{eqnarray*}
 \mbox{is a solution to}
 \begin{array}{ccc}
 &&\hskip-15mm\varphi_{,ij}-\varphi=0\ (i\ne j)
 \ \ \Leftrightarrow\ \
 \varphi=a\,e^{\sum_{j=1}^p x_j} + b\,e^{-\sum_{j=1}^p x_j},\\
 &&\hskip-6mm \varphi_{,ij}{+}\varphi=0\ (i\ne j)
 \ \Leftrightarrow\
 \varphi= a \cos(\sum\nolimits_{j=1}^p x_j) {-} b \sin(\sum\nolimits_{j=1}^p x_j),
 \end{array}
\end{eqnarray*}
where $a,b\in\RR$ and $a^2+b^2>0$.}

\vskip1mm
\textbf{Proof of Theorem~\ref{T-main-4}}.
We set $\tilde g=\bar g/\phi^2=g/(\phi F)^2=g/{\varphi}^2$, where $g$ is the Euclidean metric, and ${\varphi}=\phi F$. From (\ref{E-PDE-c1})(a), (\ref{E-Rij3a}), in view of $\ric_{i}=0$ for $g$, we~have (\ref{E-Rij3aH}). By conditions, we get the system of PDE's
\begin{eqnarray}\label{E-t4-1}
\nonumber
 p_2\,\varphi_{,x_i x_i}\eq{\varphi} f_{ii} -\Delta^{(2)}{\varphi}+{p_2}\,|\nabla{\varphi}|^2/{\varphi}
 \quad (i\le p_1),\\
 p_1\varphi_{,x_\alpha x_\alpha}\eq{\varphi}\,f_{\alpha\alpha} -\Delta^{(1)}{\varphi}+p_1|\nabla{\varphi}|^2/{\varphi}
  \quad (\alpha>p_1),\\
\nonumber
 p_2\,\varphi_{,x_i x_j}\eq f_{ij}\,\varphi\quad (1\le i\ne j\le p_1),\\
 \nonumber
 p_1\varphi_{,x_\alpha x_\beta}\eq f_{\alpha\beta}\,\varphi\quad (p_1<\alpha\ne\beta\le n),\\
 \nonumber
 \varphi_{,x_i x_\alpha}\eq 0\quad (1\le i\le p_1,\ p_1<\alpha\le n).
\end{eqnarray}
Since $T$ is non-diagonal, there is a pair $i_0,j_0$ (assume them $\le p_1$) such that $f_{i_0j_0}\ne0$ on an open set $V_1\subset V$.
From (\ref{E-t4-1})$_{3,4}$, since $f_{ij}$ are functions of two variables, and $p_1,p_2\ge3$, it follows
\begin{equation}\label{E-ff40}
  f_{ij}\varphi_{,k}= f_{ik}\varphi_{,j}=f_{jk}\varphi_{,i}\ (=p_2\varphi_{,ijk}),\quad
  \mbox{for all } i,j,k \mbox{ distinct}.
\end{equation}
\underline{If $f_{i_0 k}\equiv0$ on $V_1$, for all $k$ distinct from $i_0$ and $j_0$},
then we may assume (under a change of indices) that $f_{12}(x_1,x_2)\ne0$ and $f_{1k}\equiv0$ for $3\le k\le n$ on $V_1$.
By (\ref{E-ff40}) and (\ref{E-t4-1})$_{5}$, we have $f_{12}\varphi_{,\alpha}=p_2\varphi_{,12\alpha}=p_2(\varphi_{,1\alpha})_{,2}=0$ for $\alpha>p_1$,
moreover,
$ f_{12}\varphi_{,k}= f_{1k}\varphi_{,2}=f_{2k}\varphi_{,1}$ for $k\ge3$.
Hence, $\varphi_{,k}=f_{2k}=0$ for all $3\le k\le n$.

Observe that $\varphi_{,1}$ and $\varphi_{,2}$ cannot be zero on any open subset of $V_1$, otherwise we would have $\varphi_{,12} = f_{12}\varphi/p_2=0$. This is a contradiction since $\varphi$ is a non-vanishing function. Therefore, there exists an open subset
$V_2\subset V_1$, where $\varphi_{,1}\ne0$ and $\varphi_{,2}\ne0$.
Hence, $f_{2k}\equiv0$ on $V_2$ for all $k\ge3$. From (\ref{E-ff40}) we get
$f_{2j}\varphi_{,k} = f_{jk}\varphi_{,2}$, for $3\le j\ne k\le p_1$ and therefore $f_{jk}\equiv0$ on~$V_2$.
Differentiating (\ref{E-t4-1})$_4$, yields $f_{\alpha\beta}\varphi_{,1}=p_1 \varphi_{,\alpha\beta1}=0$, and using $\varphi_{,1}\ne0$, we conclude that $f_{\alpha\beta}\equiv0$ for all $p_1<\alpha\ne\beta\le n$ on~$V_2$.
By the above, $\varphi$ depends on $x_1$ and $x_2$ only, and is a solution to the PDE $\varphi_{,x_1x_2}=(f_{12}/p_2)\varphi$ in the $(x_1,x_2)$-plane.
Moreover, (\ref{E-t4-1})$_{1,2}$ determine the diagonal elements of $T$ which will depend on $(x_1,x_2)$ only.

\underline{Otherwise}, there exist distinct indices $i, j, k$ (assume them $\le p_1$) such that $f_{ij}$ and $f_{ik}$ do not vanish on an open subset $V_1$ of $V$. Observe that $\varphi_{,k}$ and $\varphi_{,j}$ cannot be zero on any open subset of $V_1$, since $\varphi$ is a non-vanishing differentiable function,
see~(\ref{E-t4-1})$_3$. Let $V_2\subset V_1$ be an open subset where $\varphi_{,k}\ne0$ and $\varphi_{,j}\ne 0$. It follows from (\ref{E-ff40}), $f_{jk}\ne0$ and $\varphi_{,i}\ne0$ on $V_2$. By reordering the variables, if necessary, we may consider $i = 1$ and $f_{1j}\ne0$, on an open subset $V_3\subset V_2$, for all $j$, such that $2\le j\le p$, where $p$ is an integer $3\le p\le p_1$ and $f_{1s}\equiv0$, on $V_3$ for $p+1\le s\le n$.
 Since, $\varphi$ is a non-vanishing function, there is an open subset $V_4$ of $V_3$, where
$\varphi_{,j}\ne0$ for $j = 1,\ldots, p$. It~follows from (\ref{E-ff40}) that on $V_4$,
\begin{eqnarray*}
 f_{1j}\varphi_{,k}\eq f_{jk}\varphi_{,1},\quad j\ne k,\ 2\le j,k\le p,\\
 f_{12}\varphi_{,s}\eq f_{1s}\varphi_{,2},\quad p+1\le s\le n,\\
 f_{kj}\varphi_{,s}\eq f_{sj}\varphi_{,k},\quad j\ne k,\ 2\le j,k\le p,\ p+1\le s\le n,\\
 f_{ks}\varphi_{,r}\eq f_{sr}\varphi_{,k},\quad s\ne r,\ p+1\le s,r\le n.
\end{eqnarray*}
From the first equality we get that $f_{jk}\equiv0$ on $V_4$. From the second one we
conclude that $\varphi_{,s}\equiv 0$ on $V_4$. It follows from the third one that $f_{sj}\equiv0$ and from the last equality we conclude that $f_{sr}\equiv0$ on $V_4$. Hence, $\varphi$ depends on the variables $x_1,\ldots, x_p$, and it satisfies the differential equation (\ref{E-t4-1})$_3$ for $1\le i\ne j\le p$ ,where all $f_{ij}$ do not vanish on $V_4$.
It follows from Lemma~A that, on each connected component $W\subset V_4$,
where $1\le i\ne j\le \prod\nolimits_{\,i\ne j\ne k} f_{ij}\varphi_{,k}\ne0$, there exist nonconstant differentiable functions $U_i(x_i)$, $1\le i\le p$ such that
$$
 f_{ij}/p_2 = \epsilon\,U_i'(x_i)\,U_j'(x_j),\quad 1\le i\ne j\le p_1.
$$
where $\epsilon=1$ or $\epsilon=-1$ for all $i\ne j$. We now consider on $W$ the change of variables $y_i = U_i(x_i)$. In these new coordinates $\varphi(y_1,\ldots, y_p)$ satisfies PDE's
 $\varphi_{y_i y_j} =\epsilon\,\varphi\ \mbox{for all } i\ne j$.
Lemma~B implies that $\varphi$ is given by (\ref{E-P06}) or (\ref{E-P07}) on $W$,
according to the value of $\epsilon$.
Moreover, the diagonal elements of the tensor $T$, $f_{ii}(x_1,\ldots, x_p)$ are
determined by (\ref{E-t4-1})$_1$.
In both cases, one can extend the domain of $\varphi$ to a subset of $V$ where the
functions $U_i$ are defined and $\varphi$ does not vanish.
The~converse in both cases is a straightforward computation.\hfill$\square$

\vskip1mm
\textbf{Proof of Corollary~\ref{C-main1}}.
For all cases (i)--(iv) we define $u=\varphi^{\,2/(n-2)}$.

(i) By Theorem~\ref{T-main-1}, the mixed scalar curvature for the metric $\tilde g$ is
$\widetilde K_{1,2}=-p_1p_2[\lambda+2(a_2-a_1)\,\mu]$.
Substituting  in (\ref{E-PDE-RnKK}) with $K=0$, we get (\ref{E-PDE-RnKK}).

(ii) It follows from (\ref{E-T4-fh-a})(a) that
$\widetilde K_{1,2}=p_2\,e^{2\/U}[U''-(p_1{-}1){U'}^2]$ for the metric $\tilde g$ of Theorem~\ref{T-main-2}. Similarly to (i), we obtain (\ref{E-PDE-RnKK}).

 (iii) It follows from (\ref{E-Rijabc46})(a) that
 $\widetilde K_{1,2}=(v+w)(p_2 {v}''{+}p_1{w}'') - p_1 p_2 ({v'}^{2}{+}{w'}^{2})$
 for the metric $\tilde g$ of Theorem~\ref{T-main-3}. Similarly to (i), we obtain (\ref{E-PDE-RnKK}).

(iv) Using $\varphi_{,i}=U_i'(af-bf^{-1})$ and $\varphi_{,ii}=(U_i''+{U_i'}^2)(af-bf^{-1})$,
we get $|\nabla\varphi|^2=\sum_j {U_i'}^2(af-bf^{-1})^2$
and $\Delta^{(1)}\varphi=\sum_j(U_i''+{U_i'}^2)(af-bf^{-1})$.
From (\ref{E-P06}) we get the required $\tilde K_{1,2}$ for $\tilde g$ of Theorem~\ref{T-main-4}. Hence $u$ satisfies (\ref{E-PDE-RnKK}).\hfill$\square$

\vskip1mm
Next we will prove results of Section~\ref{sec:scalar}.

\vskip1mm
Given $p_1$-dimensional plane ${\mathcal{P}}_1$ in $\RR^n$, denote by $U({\mathcal{P}}_1)$ the set of all $p_1$-dimensional planes of $\RR^n$ uniquely projecting onto ${\mathcal{P}}_1$.
Taking orthonormal basis $e_i\ (1\le i\le p_1)$ of ${\mathcal{P}}_1$, and extending it to orthonormal basis $e_i\ (1\le i\le n)$ of $\RR^n$, we represent any $\tilde{\mathcal{P}}\in U({\mathcal{P}}_1)$
as a linear graph over ${\mathcal{P}}_1$ with values in orthogonal complement, i.e.,
by the system $x_j=\sum_{i=1}^{p_1}a_{ji} x_i\ (p_1<j\le n)$.
The $p_1p_2$ elements of the matrix $A=(a_{ji})$ can be chosen as local coordinates on real Grassmannian $G_{p_1}(\RR^n)$ in a neighborhood $U({\mathcal{P}}_1)$, in particular, $\dim G_{p_1}(\RR^n)=p_1p_2$.
To any one-parameter variation ${\mathcal{P}}_1(s)$ of ${\mathcal{P}}_1$ there corresponds the matrix-function
$A(s)=(a_{ji}(s))$. Assuming $A(s)$ of a class $C^2$ with $A_1=(d A/ds)(0)$ and $A_2=(d^2 A/ds^2)(0)$, we have $A(s)=sA_1+\frac12 s^2 A_2+o(s^2)$.
For $\tilde{\mathcal{P}}\in U({\mathcal{P}}_1)$, the stationary values $0\le\alpha_1\le\ldots\le\alpha_{p_1}\le\pi/2$
of angle between unit vectors in $\tilde{\mathcal{P}}$ and ${\mathcal{P}}_1$ are called the \textit{angles} between these planes. It is known that there exist orthonormal bases
$e_i\ (1\le i\le p_1)$ of ${\mathcal{P}}_1$ and $\tilde e_j\ (1\le j\le p_1)$ of $\tilde{\mathcal{P}}$
such that $\langle e_i, \tilde e_j\rangle=\cos\alpha_i\delta_{ij}$
(the property can be taken as the definition of angles $\alpha_i$).
The angles $\alpha_i\ (1\le i\le p_1)$ determine the relative position of two  $p_1$-dimensional planes in $\RR^n$,
and $\rho(\tilde{\mathcal{P}}, {\mathcal{P}}_1)=(\sum_{i=1}^{p_1}\alpha_i^2)^{1/2}$ is the \textit{distance}.
Let $p_1\le p_2$.
Choosing the basis, one may represent a variation ${\mathcal{P}}_1(s)$ of ${\mathcal{P}}_1$
using the matrix $A(s)=sA_1+\frac12 s^2 A_2+o(S^2)$, where $A_1$ consists of
$a^{(1)}_{ji}=\delta_{ji}\cos\alpha_i\ (1\le i\le p_1,\ p_1<j\le n)$.

\textbf{Proof of Theorem~\ref{T-Kmix}}.
Let $\varphi:M\to\RR$ be a smooth nonnegative function with a local support $V\subset M$.
Assume that $p_1\le p_2$. Due to discussion above, any variation ${\mathcal{D}}_1(s)\ (|s|<1)$ on $V$ can be represented by orthogonal vector fields
$e_{i}(s)=e_{i}\cos(s\,\omega_i\,\varphi)+\eps_i\sin(s\,\omega_i\,\varphi)$,
where
$e_i\ (i\le p_1)$ is orthonormal basis of ${\mathcal{D}}_{1|\,V}$,
$\eps_j\ (j\le p_2)$ orthonormal basis in ${\mathcal{D}}_2$ on $V$,
and $\omega_i\ge0\ (i\le p_1)$ are real numbers.
 Set $\eps_{i}(s)=\eps_{i}\cos(s\,\omega_i\,\varphi)-e_i\sin(s\,\omega_i\,\varphi)$
for $1\le i\le p_1$ and $\eps_j(s)=\eps_j$ for $p_1< j\le n$.
Indeed, $\eps_j(s)\ (1\le j\le p_2)$ span ${\mathcal{D}}_2(s)$ on~$V$.
 We have $\frac{de_i}{ds}(0)=\varphi\,\omega_i\,\eps_i$ and
$\frac{d\eps_i}{ds}(0)=-\varphi\,\omega_i\,e_i$ on~$V$ for $i\le p_1$.

We extend distributions outside of $V$ as ${\mathcal{D}}_1(s)={\mathcal{D}}_1$ and define $I(s)=I_K({\mathcal{D}}_1(s))$ for all $s$. The mixed scalar curvature of ${\mathcal{D}}_1(s)$ over $V$ is
$$
 K(s):=K({\mathcal{D}}_1(s),{\mathcal{D}}_2(s))= \sum\nolimits_{i=1}^{p_1}\sum\nolimits_{j=1}^{p_2}
 g(R(e_i(s), \eps_j(s))\eps_j(s), e_i(s)).
$$
Derivation by $s$ at $s=0$ (using symmetries of curvature tensor) yields on $V$
\begin{eqnarray*}
 \frac{d}{ds}K(s)_{|\,s=0}\eq 2\,\varphi\sum\nolimits_{i,j=1}^{p_1}
 [g(R(e_i, \eps_j)\eps_j, \omega_i\eps_i)+g(R(e_i, \eps_j)(-\omega_j e_j), e_i)]\\
 \eq 2\,\varphi\sum\nolimits_{i=1}^{p_1}\omega_i(\ric_1-\ric_2)(e_i, \eps_i).
\end{eqnarray*}
Notice that $K(s)=K(0)$ outside of $V$. Hence
\begin{equation*}
 I_K'(0) = 2\int_V\varphi\sum\nolimits_{i=1}^{p_1}\omega_i(\ric_1-\ric_2)(e_i, \eps_i)\,{\rm d}\vol.
\end{equation*}
Since $e_i,\eps_i,\,\omega_i$ and $\varphi$ are arbitrary,
from above it follows that ${\mathcal{D}}$ is critical for $I_K$ if and only if
$(\ric_1-\ric_2)(X, Y)=0$ for all $X\in{\mathcal{D}_1}$ and $y\in{\mathcal{D}_2}$, the first part of the claim has been proved.

Assume now that ${\mathcal{D}}_1$ is a critical point for $I_K$.
Consider arbitrary variation ${\mathcal{D}}_1(s)\ (|s|<1)$ of ${\mathcal{D}}_1$ with a local support $V\subset M$.
As above, take orthonormal basis $e_{i}\ (i\le p_1)$ of ${\mathcal{D}}_1$ and orthonormal basis $\eps_i\ (i\le p_2)$ of ${\mathcal{D}}_2$ on $V$ such that ${\mathcal{D}}_1(s)$ is spanned by $e_{i}(s)$ and ${\mathcal{D}}_2(s)$ is spanned by $\eps_{j}(s)$ as above, where $\varphi:M\to\RR$ a smooth nonnegative function with a support in~$V$.
In addition, $\frac{d^2e_i}{ds^2}(0)=-\varphi^2\,\omega_i^2\, e_i$ and
$\frac{d^2\eps_i}{ds^2}(0)=-\varphi^2\,\omega_i^2\,\eps_i$ on~$V$  for $i\le p_1$.
The second derivative of $K(s)$ at $s=0$ on $V$ ($K(s)=K(0)$ outside of $V$) is
\begin{eqnarray*}
 \frac12\,
 K''_{|\,s=0} 
 \eq\sum\nolimits_{i,j=1}^{p_1}
 \big[\,g(R(e_i,\eps_j)\eps_j, e_i'')+g(R(e_i,\eps_j)\eps_j'', e_i)+g(R(e_i',\eps_j)\eps_j, e_i')\\
 \plus g(R(e_i,\eps_j')\eps_j', e_i)+2\,g(R(e_i,\eps_j')\eps_j, e_i')+2\,g(R(e_i,\eps_j)\eps_j', e_i')\,\big]\\
 \eq\varphi^2\big\{
  \sum\nolimits_{j=1}^{p_1}\omega_j^2(\ric_{1}-\ric_{2})(\eps_j,\eps_j)
 -\sum\nolimits_{i=1}^{p_1}\omega_i^2(\ric_{1}-\ric_{2})(e_i,e_i)\\
 \plus2\sum\nolimits_{i,j=1}^{p_1}
 \omega_i\,\omega_j[\,g(R(e_i,e_j)\eps_i, \eps_j)+g(R(e_i,\eps_j)\eps_i, e_j)\,]\big\}.
\end{eqnarray*}
Hence
 $I_K''(0) = 2\int_V\varphi^2\sum\nolimits_{i,j=1}^{p_1}
  \Phi_{ij}(x)\,\omega_i\,\omega_j\,{\rm d}\vol$,
where $\Phi_{ij}(x)$ is defined in (\ref{E-phi-ij}) and $\omega_i\in\RR$ are arbitrary numbers. From above it follows that $I_K''(0)>0$ when the form ${\cal I}$ is quasi-positive.$\hfill\square$

\vskip1mm
\textbf{Proof of Proposition~\ref{L-IntK}}.
Define the linear operators $C_{2}^{\,i}:{\mathcal{D}}_2\to {\mathcal{D}}_2$ by
$C_2(N,\cdot)=\sum\nolimits_{\,i}N_i C_{2}^{\,i}$
for $N=\sum\nolimits_{\,i} N_i e_i\in {\mathcal{D}}_1$. Note that
$p_2H_2=-\sum\nolimits_i\sigma_1(C_{2}^{\,i})\,e_i$. Hence
 $p_2^{\,2}|H_2|^2=\sum\nolimits_i \sigma_1^2(C_{2}^{\,i})$
 and
 $\sigma_1((C_{2}^{\,i})^2)=\sum\nolimits_{\alpha,\beta}
 C_{2,\alpha\beta}^{\,i}C_{2,\beta\alpha}^{\,i}$.
\newline
From
 $-g(\nabla_{e_\alpha}e_\beta, e_i)=g(e_\beta, \nabla_{e_\alpha}e_i)=C_{2,\alpha\beta}^{\,i}$,
we get
$(\nabla_{e_\alpha}e_\beta)^\top{=}-\sum_{i}C_{2,\alpha\beta}^{\,i}\,e_i$.
Hence
 \begin{equation*}
 \begin{array}{c}
 |B_2|^2-|T_2|^2=\sum\nolimits_{\alpha,\beta} g((\nabla_{e_\alpha}e_\beta)^\top, (\nabla_{e_\beta}e_\alpha)^\top)=\\
 \sum\nolimits_{\alpha,\beta} g\big(\sum\nolimits_{i}
 C_{2,\alpha\beta}^{\,i} \,e_i, \sum\nolimits_{j} C_{2,\alpha\beta}^{\,j} \,e_j\big)
 =\sum\nolimits_{i}\sum\nolimits_{\alpha,\beta} C_{2,\alpha\beta}^{\,i} C_{2,\beta\alpha}^{\,i}=\sum\nolimits_{i}\sigma_1((C_{2}^{\,i})^2).
 \end{array}
 \end{equation*}
 From known identity $2\,\sigma_2(C_{2}^{\,i})=\sigma_1^2(C_{2}^{\,i})-\sigma_1((C_{2}^{\,i})^2)$,
 it follows from above that
 $2\sum\nolimits_i \sigma_2(C_{2}^{\,i})=p_2^{\,2}|H_2|^2+|T_2|^2-|B_2|^2$.
 Similarly, for $Y=\sum\nolimits_{\,\alpha} Y_\alpha e_\alpha$ we get $C_1(Y,\cdot)=\sum\nolimits_{\,\alpha} Y_\alpha C_{1}^{\,\alpha}$.
 Thus,
 $2\sum\nolimits_\alpha\sigma_2(C_{1}^{\,\alpha})=p_1^{\,2}|H_1|^2+|T_1|^2-|B_1|^2$
 and
\begin{equation}\label{E-2parts1}
  2\hskip-2pt\sum\nolimits_i\sigma_2(C_{2}^{\,i})
 +2\hskip-2pt\sum\nolimits_\alpha\sigma_2(C_{1}^{\,\alpha})
 = p_1^{\,2}|H_1|^2 +|T_1|^2-|B_1|^2+p_2^{\,2}|H_2|^2+|T_2|^2-|B_2|^2.
\end{equation}
Using Lemma~1.5 of \cite{rw3}, we obtain at $x\in M$
\begin{equation}\label{E-3parts}
 \sigma_2({\mathcal{D}}_1)+\sigma_2({\mathcal{D}}_2)=
 \sum\nolimits_{\,i}\sigma_2(C_{2}^{\,i})+\sum\nolimits_{\,\alpha}\sigma_2(C_{1}^{\,\alpha}).
\end{equation}
The integral formula in \cite{w90} shows us that
\begin{equation}\label{eq-wal2}
\begin{array}{c}
 \int_{M}\big[K_{1,2}+|B_1|^2{-}p_1^{\,2}|H_1|^2{-}|T_1|^2{+}|B_2|^2{-}p_2^{\,2}|H_2|^2{-}|T_2|^2\big]
 \,{\rm d}\vol=0,
\end{array}
 \end{equation}
where
$H_i,\,T_i,\,B_i$ are the mean curvature vector, the integrability tensor and the 2-nd fundamental form of ${\mathcal{D}}_i$. The required formula
follows directly from (\ref{E-2parts1})\,--\,(\ref{eq-wal2}).$\hfill\square$

\vskip1mm
\textbf{Proof of Proposition~\ref{P-bend}}.
We represent the norm of $\nabla\tilde {\mathcal{D}}_1$ using co-nullity operators $C_1,\,C_2$
\begin{equation}\label{E-bendD0}
 |\nabla\tilde {\mathcal{D}}_1|^2=\sum\nolimits_{\,i,j,\alpha}(C_{1,ij}^{\,\alpha})^2
 +\sum\nolimits_{\,i,\alpha,\beta}(C_{2,\alpha\beta}^{\,i})^2.
\end{equation}
Hence $|\nabla\tilde {\mathcal{D}}_1|^2=|\nabla\tilde {\mathcal{D}}_2|^2$.
Clearly, this expression does not depend on the adapted local orthonormal basis.
Let $\tilde C=\{\tilde C_{ij}\}_{1\le i,j\le p}$ be a
real matrix of order $p\ge2$. An elementary calculation shows that
\begin{equation}\label{E-bendD2}
\begin{array}{c}
 (p-1)\sum_{\,i,j}(\tilde C_{ij})^2=
 \sum_{\,i<j}(\tilde C_{ii}-\tilde C_{jj})^2+\sum_{\,i<j}(\tilde C_{ij}+\tilde C_{ji})^2\\
 +(p-2)\sum_{\,i\ne j}(\tilde C_{ij})^2
 +2\sum_{\,i<j}(\tilde C_{i\,i}\,\tilde C_{jj}-\tilde C_{ij}\,\tilde C_{ji}).
\end{array}
\end{equation}
Note that the last term in (\ref{E-bendD2}) equals $\sigma_2(\tilde C)$. Hence
 $\sum\nolimits_{\,i,j}(\tilde C_{ij})^2\ge \frac2{p-1}\,\sigma_2(\tilde C)$.
Applying this inequality to the matrices $C_1^{\,\alpha}$ and $C_2^{\,i}$, we obtain
\begin{equation*}
\begin{array}{c}
 \sum\nolimits_{\,\alpha;\,i,j}(C_{1,ij}^{\,\alpha})^2\ge \frac2{p_1{-}1}\sum\nolimits_{\,\alpha}\sigma_2(C_{1}^{\,\alpha}),\quad \sum\nolimits_{\,i;\,\alpha\,\beta}(C_{2,\alpha\,\beta}^{\,i})^2\ge
 \frac2{p_2{-}1}\sum\nolimits_{\,i}\sigma_2(C_{2}^{\,i}),
 \end{array}
\end{equation*}
and, in view of (\ref{E-bendD0}),
 $|\nabla\tilde {\mathcal{D}}_1|^2\ge
 \frac2{p_1{-}1}\sum\nolimits_{\,\alpha}\sigma_2(C_{1}^{\,\alpha})
 +\frac2{p_2{-}1}\sum\nolimits_{\,i}\sigma_2(C_{2}^{\,i})$.
From that the inequality (\ref{E-bendD5}) follows.$\hfill\square$


\end{document}